\documentclass[12pt]{article}
\usepackage{amsmath,amsfonts,amscd,amsxtra}

\begin{document}

\newcommand{\proof}{{\bf \underline{\underline{Proof}} : }}

\newcommand{\tit}[1]{\bigskip\begin{bf} \begin{center} \begin{Large}
\arabic{section}. #1
\end{Large}\end{center}\end{bf}}

\newcommand{\bigtit}{\begin{bf} \begin{center} \begin{LARGE}
``Torsion and fibrations''  \\
by\\
Wolfgang L\"uck, Thomas Schick and Thomas Thielmann
\end{LARGE}\end{center}\end{bf}}

\newcommand{\cat}{{\cal C}}

\newcommand{\F}{{\cal F}}

\newcommand{\qed}{\hspace{10mm} \rule{2mm}{2mm}}

\newcommand{\tr}{\operatorname{tr}}
\newcommand{\pr}{\operatorname{pr}}
\newcommand{\id}{\operatorname{id}}
\newcommand{\im}{\operatorname{im}}
\newcommand{\cok}{\operatorname{cok}}
\newcommand{\vol}{\operatorname{vol}}
\newcommand{\clos}{\operatorname{clos}}
\newcommand{\cone}{\operatorname{cone}}
\newcommand{\cyl}{\operatorname{cyl}}
\newcommand{\aut}{\operatorname{aut}}
\newcommand{\Wh}{\operatorname{Wh}}
\newcommand{\sing}{\operatorname{sing}}
\newcommand{\cell}{\operatorname{cell}}
\newcommand{\SPACES}{\operatorname{SPACES}}
\newcommand{\COMPLEXES}{\operatorname{COMPLEXES}}
\newcommand{\COCHAIN}{\operatorname{COCHAIN}}
\newcommand{\CHAIN}{\operatorname{CHAIN}}
\newcommand{\VECTOR}{\operatorname{VECTOR}}
\newcommand{\trans}{\operatorname{trans}}
\newcommand{\odd}{\operatorname{odd}}
\newcommand{\ev}{\operatorname{ev}}
\newcommand{\fil}{\operatorname{fil}}
\newcommand{\Serre}{\operatorname{LS}}
\newcommand{\dR}{\operatorname{dR}}
\newcommand{\an}{\operatorname{an}}
\newcommand{\Pf}{\operatorname{Pf}}
\newcommand{\topological}{\operatorname{top}}
\newcommand{\RS}{\operatorname{RS}}
\newcommand{\harm}{\operatorname{harm}}

\newcommand{\rr}{\mathbb R}
\newcommand{\zz}{\mathbb Z}

\newcommand{\brn}{\begin{refnumber} \em }
\newcommand{\ern}{\em \end{refnumber}}

\newtheorem{theorem}{Theorem}[section]
\newtheorem{proposition}[theorem]{Proposition}
\newtheorem{lemma}[theorem]{Lemma}
\newtheorem{definition}[theorem]{Definition}
\newtheorem{example}[theorem]{Example}
\newtheorem{diagram}[theorem]{Diagram}
\newtheorem{remark}[theorem]{Remark}
\newtheorem{notation}[theorem]{Notation}
\newtheorem{refnumber}[theorem]{}
\newtheorem{corollary}[theorem]{Corollary}
\newtheorem{assumption}[theorem]{Assumption}
\newtheorem{conjecture}[theorem]{Conjecture}
\newtheorem{problem}[theorem]{Problem}

{\catcode`@=11\global\let\c@equation=\c@theorem}
\renewcommand{\theequation}{\thetheorem}


\newcommand{\comsquare}[8]{
\begin{center}
$\begin{CD}
#1 @>#2>> #3\\
@V{#4}VV @VV{#5}V\\
#6 @>>#7> #8
\end{CD}$
\end{center}}

\newcommand{\shortexactsequence}[5]
{0 \xrightarrow{} #1 \xrightarrow{#2} #3 \xrightarrow{#4} #5
\xrightarrow{}  0}

\newcommand{\squarematrix}[4]{\left( \begin{array}{cc} #1 & #2 \\ #3 & #4
\end{array} \right)}



\typeout{-------------------------- Introduction
 --------------------------}
\bigtit

\begin{abstract}
We study the behaviour of analytic torsion
under smooth fibrations. Namely, let
\mbox{$F \rightarrow E \xrightarrow{f} B$}
be a smooth fiber bundle of connected
closed oriented smooth manifolds and let $V$ be a flat vector bundle
over $E$. Assume that $E$ and $B$  come with Riemannian metrics.
Suppose that $\dim(E)$ is odd and $V$ is unimodular and comes with an
arbitrary Riemannian metric or that $\dim(E)$ is even and $V$ comes
with a unimodular (not necessarily flat) Riemannian metric.
Let $\rho_{\an}(E;V)$ be the analytic torsion of $E$ with coefficients
in $V$, let $\rho_{\an}(F_b;V)$ be the analytic torsion of the fiber
over $b$ with coefficients in $V$ restricted to $F_b$ and
let $\Pf_B$ be the Pfaffian $\dim(B)$-form. Let $H^q_{\dR}(F;V)$
be the flat vector bundle over $B$ whose fiber over $b \in B$ is
$H^q_{\dR}(F_b;V)$ with the Riemannian metric which comes from
the Hodge-deRham decomposition and the Hilbert space structure
on the space of harmonic forms induced by the Riemannian metrics.
Let $\rho_{\an}(B;H^q_{\dR}(F;V))$ be the analytic torsion of $B$
with coefficients in this bundle. The Leray-Serre spectral sequence
for deRham cohomology determines a certain correction term
$\rho^{\Serre}_{\dR}(f)$.
We prove
$$\rho_{\an}(E;V) ~ = ~ \int_B \rho_{\an}(F_b;V) \cdot \Pf_B
 ~ + ~ \sum_{q} (-1)^q \cdot \rho_{\an}(B;H^q_{\dR}(F;V))
 ~ + \rho^{\Serre}_{\dR}(f).$$
This formula simplifies in special cases such as bundles
with $S^n$ as fiber or base, in which case the correction terms
$\rho^{\Serre}_{\dR}(f)$ reduces to the torsion of the
associated Gysin or  Wang sequence, resp.
\end{abstract}
\tit{Introduction}

Let $M$ be a connected closed smooth manifold with Riemannian metric
and $V$ be a flat vector bundle with a not necessarily flat Riemannian
metric. The definition of analytic torsion
due to  Ray and Singer \cite{Ray-Singer (1971)} for an orthogonal
representation $V$, or,
equivalently, for a flat Riemannan metric on $V$,
still makes sense in the setting above
\cite[page 35]{Bismut-Zhang (1992)}
and \cite[page 730]{Mueller 1993}. Namely,
let $\zeta_p(s)$ be the zeta-function of the Laplace operator
\mbox{$\Delta^p : \Omega^p(E;V) \longrightarrow \Omega^p(E;V)$}
which is for $\mbox{Re}(s) >> 0$ the holomorphic function
\mbox{$\sum_{\lambda > 0} \lambda^{-s}$} where $\lambda$ runs over the
positive eigenvalues of $\Delta^p$. It has a meromorphic extension to
the complex plane which is analytic in zero. Define
\begin{eqnarray}
\rho_{\an}(E;V) & := & \frac{1}{2} \cdot
\sum_{q \ge 0} (-1)^q \cdot q \cdot \zeta_p^{\prime}(0)
\hspace{5mm}\in \rr.
\label{analytic torsion}
\end{eqnarray}
We want to study it for
smooth fiber bundles. The main result of this paper is
\begin{theorem} \label{fiber bundle formula for analytic torsion}
\hspace*{-1pt}
Let \mbox{$F \rightarrow E \xrightarrow{f} B$}
be a smooth fiber bundle of connected
closed oriented smooth manifolds and let $V$ be a flat vector bundle
over $E$. Assume that $E$ and $B$  come with Riemannian metrics.
Suppose that $\dim(E)$ is odd and $V$ is unimodular and comes with an
arbitrary Riemannian metric or that $\dim(E)$ is even and $V$ comes
with a unimodular Riemannian metric. Then
$$\rho_{\an}(E;V) ~ = ~ \int_B \rho_{\an}(F_b;V) \cdot \Pf_B
 ~ + ~ \sum_{q} (-1)^q \cdot \rho_{\an}(B;H^q_{\dR}(F;V))
 ~ + \rho^{\Serre}_{\dR}(f). \qed $$
\end{theorem}

Here are some explanations of the assumptions and the formula in
Theorem \ref{fiber bundle formula for analytic torsion}. \par

For a path $w$  in $E$ the fiber transport gives a linear
isomorphism \mbox{$V_w : V_{w(0)} \longrightarrow V_{w(1)}$}
which depends only on the homotopy class relative endpoints of $w$
since $V$ is flat. We call $V$ {\em unimodular} if for one (and hence all) $e \in E$ and all loops $w$ with base point $e$ we get
\mbox{$| \det(V_w : V_e \longrightarrow V_e)| = 1$}.
We call a Riemannian metric on $V$
{\em unimodular} if for any path $w$ in $E$ we get
\mbox{$\det(V_w^{\ast}\circ V_w : V_{w(0)}
\longrightarrow V_{w(0)})   =  1$}
where $V_{w}^{\ast}$ is the adjoint of $V_w$ with respect to the
Hilbert space structure on the
fibers of $V$ given by the Riemannian metric.
This is a weaker condition than being a flat Riemannian metric
what would mean that $V_w$ is always an isometry.
Notice that $V$ is
unimodular if and only if it carries a unimodular
Riemannian metric. \par

Of course $\rho_{\an}(E;V)$ is just the analytic torsion with respect to
the given Riemannian metrics on $E$ and $V$. 
These induce also a metric on $H^p_{\dR}(E,V)$.
Each fiber $F_b = p^{-1}(b)$ inherits a Riemannian metric
from $E$ by restriction. We denote the restriction of $V$ to $F_b$
again by $V$. Hence $\rho_{\an}(F_b;V)$ is defined
and is a smooth function in $b \in B$.  \par

Let $\Pf_B$ be the Pfaffian $\dim(B)$-form on the oriented
Riemannian manifold $B$. It is a representative of the Euler
class of $B$ in Chern-Weil theory and satisfies by the Gauss-Bonnet
theorem
$$\int_B \Pf_B ~ = ~ \chi(B)$$
where $\chi(B)$ is the Euler characteristic. If $\dim(B)$  is odd,
then $\Pf_B$ is defined to be zero. \par

Let ${\cal H}^q(F_b;V)$ be the space of harmonic $q$-forms,
i.e. the kernel of the Laplace operator
\mbox{$\Delta^q : \Omega^q(F_b;V) \longrightarrow \Omega^q(F_b;V)$}.
It inherits an inner product from the Riemannian
metrics on $F_b$ and $V$. The {\em harmonic Hilbert structure} on the
deRham cohomology $H^q_{\dR}(F_b;V)$ is the Hilbert space structure
for which the canonical Hodge isomorphism
$${\cal H}^q(F_b;V) \longrightarrow H^q_{\dR}(F_b;V)$$
is isometric. Thus we get a Riemannian metric on the canonical flat
vector bundle $H^q_{\dR}(F;V)$ whose fiber over
$b \in B$ is $H^q_{\dR}(F_b;V)$. Hence the analytic torsion
\mbox{$\rho_{\an}(B;H^q_{\dR}(F;V))$} is defined,
and $H^p_{\dR}(B,H^q_{\dR}(F,V))$ inherits a natural Hilbert space structure.\par

There is the following natural descending filtration
of the deRham complex \mbox{$\Omega^{\ast}(E;V)$}. Define
\mbox{$F_p\Omega^{n}(E;V)$} to be those $n$-forms
with coefficients in $V$
which can be written as finite sums of $n$-forms on $E$
with coefficients in $V$ of the shape
\mbox{$\omega \wedge f^{\ast}\eta$}
for \mbox{$\omega \in \Omega^{n-k}(E;V)$}
and \mbox{$\eta \in \Omega^k(B)$} for some $k \ge p$.
This filtration is compatible with the differential since
\mbox{$d(\omega \wedge f^{\ast}\eta) =
d(\omega) \wedge f^{\ast}\eta \pm \omega \wedge f^{\ast}d(\eta).$}
The associated spectral cohomology sequence is the Leray-Serre
spectral sequence for deRham cohomology,  which we recall in
Section \ref{Torsion and spectral sequences}. 

 Part of
the Leray-Serre spectral sequence for deRham cohomology is
the filtration of the cohomology $H^n(E;V)$
$$\{0\} = F^{n+1,-1} \subset \ldots F^{p+1,n-p-1} \subset F^{p,n-p}
\subset \ldots \subset F^{0,n} = H^n_{\dR}(E;V),$$
the natural identification of the $E_2$-term
\begin{eqnarray}
V_2^{p,q}: H^p_{\dR}(B;H^q_{\dR}(F;V)) & \xrightarrow{\cong} &
E^2_{p,q},
\label{identification of E_2term for deRham}
\end{eqnarray}
the identification of the cohomology of the
$r$-th term of the spectral sequence
with the $(r+1)$-th term and the identification of the
$E_{\infty}$-term with the filtration quotients
$$\begin{array}{llll}
\phi_r^{p,q}: & H^{0}(E_r^{p+r\ast,q-(r-1)\ast})
& \xrightarrow{\cong} &
 E_{r+1}^{p,q},
\\
\psi^{p,q}: & F^{p,q}/F^{p+1,q-1} & \xrightarrow{\cong} &
E_{\infty}^{p,q}.
\end{array}$$
For $r$ sufficiently large, the differentials in $E^{*,*}_r$ are trivial. 

Next we explain the
term \mbox{$ \rho^{\Serre}_{\dR}(f)$} appearing in Theorem
\ref{fiber bundle formula for analytic torsion}.\par


For a linear isomorphism $f:V\to W$ of finite-dimensional real Hilbert spaces, set
\begin{eqnarray}
[[f]] & := & \frac{1}{2}\ln(|\det(f^*f)|) ~ \in \rr
\label{[[]]}
\end{eqnarray}

Let $C = C^{\ast}$ be an acyclic finite Hilbert cochain complex. Define
\begin{eqnarray}
\rho(C) & := &
[[(c^{\ast} + \gamma^{\ast}) : C^{\ev} \longrightarrow
C^{\odd}]] ~ \in \rr
\label{Reidemeistert torsion for acyclic chain complexes}
\end{eqnarray}
where $c^{\ast}$ is the differential and $\gamma^{\ast}$ a
chain contraction. 
If \mbox{$f : C \longrightarrow D$} is a chain homotopy
equivalence of finite Hilbert cochain complexes, 
$\cone (f)$
is the cochain complex with $n$-th differential
$$\squarematrix{c^{n}}{0}{f^n}{-d^{n-1}} : C^{n} \oplus D^{n-1}
\longrightarrow
C^{n+1} \oplus D^{n}.$$
It is acyclic and we define
\begin{eqnarray}
t(f) & := & \rho(\cone(f)).
\label{torsion of chain homotopy equivalence}
\end{eqnarray}
Let $C$ be a finite Hilbert cochain complex such that
$H(C^{\ast})$ carries a Hilbert structure.
There is up to homotopy precisely one chain map
\mbox{$i : H(C) \longrightarrow C$} whith $H(i)=id$,
 where we consider $H(C)$
as a cochain complex with the trivial differential. Define
\begin{eqnarray}
\rho(C) & ~:= & -t(i) ~ \in \rr.
\label{Milnor torsion of a chain complex}
\end{eqnarray}

In the Leray-Serre spectral sequence, equip $E_2^{p,q}$ with the Hilbert space
structure which makes $V_2^{p,q}$ to an isometry. Equip inductively $E^{p+r\ast,q-(r-1)\ast}_r$, $H(E^{p+r\ast,q-(r-1)\ast}_r)$ ($r\ge 2$) and $E_\infty^{p,q}$ with the Hilbert sub- and quotient structures. In particular, $\phi^{p,q}_r$ become isometries. Equip $F^{p,q}\subset H^{p+q}_{\dR}(E,V)$ and $F^{p,q}/F^{p+1,q-1}$ with the Hilbert sub- and quotient structures
Now define
\[  \rho^{\Serre}_{\dR}(f) := 
\sum_{r \ge 2} \sum_{p=0}^{r-1} \sum_{q } (-1)^{p+q}
\cdot \rho(E_r^{p+r\ast,q-(r-1)\ast})
~ - ~ \sum_{p,q} (-1)^{p+q} \cdot [[\psi^{p,q}]]
\]
This number depends only on $f$ and the Riemannian metrics on $E$, $B$ and~$V$.

In general the correction term
$\rho^{\Serre}_{\dR}(f)$ is very involved
and is as complicated as the Leray-Serre spectral
cohomology sequence is. However, there
are cases where $\rho^{\Serre}_{\dR}(f)$ and the whole formula in
Theorem \ref{fiber bundle formula for analytic torsion} are easy to
 understand. Namely, we will prove under the assumptions of
Theorem \ref{fiber bundle formula for analytic torsion} the following
three corollaries. The first one generalizes a result of
Fried \cite{Fried (1987)} for orthogonal $V$.


\begin{corollary} \label{acyclic case}
Suppose that $H^q_{\dR}(F;V)$ vanishes for all $q$.
Then $\chi(B) \cdot \rho_{\an}(F_b;V)$ is independent of $b$ and
$$\rho_{\an}(E;V) ~ = ~ \chi(B) \cdot \rho_{\an}(F_b;V). \qed$$
\end{corollary}

\begin{corollary} \label{spheres as fibers}
Suppose that $F$ is $S^n$ and $V = f^{\ast}W$ for a flat vector
bundle with
Riemannian metric over $B$. Let $G^{\ast}$ be the acyclic cochain
complex of finite-dimensional Hilbert spaces given by the Gysin
sequence
$$ \ldots \xrightarrow{\int} H^p_{\dR}(B;W) \xrightarrow{\wedge}
H^{p+n+1}_{\dR}(B;W)
\xrightarrow{f^{\ast}} H^{p+n+1}_{\dR}(E;V) \xrightarrow{\int}
H^{p+1}_{\dR}(B;W) \xrightarrow{\wedge e(f)} \ldots$$
where $\wedge e(f)$ is the product with the Euler
class $e(B) \in H^{n+1}_{\dR}(B)$, $\int$ is
integration over the fiber and
\mbox{$G^1 = H^0_{\dR}(E;V)$}. Then the torsion
\mbox{$\rho(G^{\ast}) \in \rr$} is defined and we get:
$$\rho_{\an}(E;V) ~ = ~
\chi(S^n) \cdot \rho_{\an}(B;W) + \rho(G^{\ast}).
\qed $$
\end{corollary}
The condition \mbox{$V = f^{\ast}W$} is no loss of generality,
provided
that $n \ge 2$ or that $f$ induces an isomorphism
\mbox{$\pi_1(E) \longrightarrow \pi_1(B)$}.

\begin{corollary} \label{spheres as base}
Suppose that $B = S^n$. Let $W^{\ast}$ be the acyclic cochain complex
of finite-dimensional Hilbert spaces given by the Wang sequence and
the harmonic structures for some $b \in B$:
$$ \ldots \longrightarrow H^{q-1}_{\dR}(E;V) \longrightarrow
H^{q-1}_{\dR}(F_b;V) \longrightarrow H^{q-n}_{\dR}(F_b;V)
\longrightarrow H^q_{\dR}(E;V) \longrightarrow \ldots$$
where $W^1 = H^0_{\dR}(E;V)$. Then
$$\rho_{\an}(E;V) ~ = ~
\chi(S^n) \cdot \rho_{\an}(F_b;V) + \rho(W^{\ast}).
\qed
$$
\end{corollary}

We make some remarks about the proof of Theorem
\ref{fiber bundle formula for analytic torsion}.
It will depend on the
following deep results of Bismut-Zhang
\cite{Bismut-Zhang (1992)} and
M\"uller \cite{Mueller 1993}.
In the sequel $M$ is a connected closed
oriented Riemannian manifold and $V$ is a flat vector
bundle over $M$ with Riemannian metric. We have introduced the
analytic torsion $\rho_{\an}{}(M;V)$ above. Its topological
counterpart
\begin{eqnarray}
\rho_{\topological}(M;V) & \in & \rr.
\label{topological torsion}
\end{eqnarray}
is the Milnor torsion of $M$ with respect to some triangulation
and the harmonic Hilbert structure on cohomology which we will
recall in Definition \ref{Milnor torsion}.

\begin{theorem} \label{equality of analytic and topological torsion}
If the Riemanniann metric on $V$ is unimodular, then
$$\rho_{\an}(M;V) ~ = ~ \rho_{\topological}(M;V). \qed $$
\end{theorem}

Let $g_M$ and $\overline{g_M}$ be two Riemannian metrics on
$M$ and $g_V$ and $\overline{g_V}$ be two arbitrary Riemannian
metrics on $V$. Let $\rho_{\an}(M;V)$ and
$\overline{\rho_{\an}(M;V)}$
be the analytic torsion with respect to \mbox{$(g_M,g_V)$} and
\mbox{$(\overline{g_M},\overline{g_V})$}. Analogously we
denote the Hilbert spaces $H^p_{\dR}(M;V)$ and
$\overline{H^p_{\dR}(M;V)}$ equipped with the harmonic Hilbert
structures with respect to \mbox{$(g_M,g_V)$} and
\mbox{$(\overline{g_M},\overline{g_V})$} and the Hilbert spaces
$V_x$ and $\overline{V_x}$
equipped with the Hilbert structure with respect to
$g_V$ and $\overline{g_V}$. Denote by $\Pf_M$ the Pfaffian
with respect to $g_M$. Let $\widetilde{\Pf}(M,g_M,\overline{g_M})$
be the Chern-Simons $n-1$-form. Its image under the differential
is the difference of the two Pfaffians of $M$ with respect to
$g_M$ and $\overline{g_M}$. Let $\theta(V,\overline{g_V})$
be the closed $1$-from defined in
\cite[Definition 4.5]{Bismut-Zhang (1992)}.
It measures the deviation of $\overline{g_V}$ from being
unimodular and vanishes if $\overline{g_V}$ is unimodular.

\begin{theorem} \label{variation of Riemanninan metrics} We get
under the conditions and in the notations above:
\begin{enumerate}
\item If $\dim(M)$ is odd, then
$$\rho_{\an}(M;V) -
\overline{\rho_{\an}(M;V)} ~ = ~ -
\sum_{p} (-1)^p \cdot [[H^p_{\dR}(M;V)
\xrightarrow{\id} \overline{H^p_{\dR}(M;V)}]];$$
\item If $\dim(M)$ is even, then
\begin{eqnarray*}
\rho_{\an}(M;V)  -
\overline{\rho_{\an}(M;V)} & = &
 -  \sum_{p} (-1)^p \cdot [[H^p_{\dR}(M;V)
\xrightarrow{\id} \overline{H^p_{\dR}(M;V)}]] \\
& &  +
\int_M [[V_x \xrightarrow{\id} \overline{V_x}]] \cdot \Pf_M  ~ - ~
\int_{M} \theta(V,\overline{g_V})
\cdot \widetilde{\Pf}(M,g_M,\overline{g_M}).
\hspace{2mm} \rule{2mm}{2mm}
\end{eqnarray*}
\end{enumerate}
\end{theorem}

Theorem \ref{equality of analytic and topological torsion}
and Theorem \ref{variation of Riemanninan metrics}
for odd-dimensional $M$ are due to M\"uller
\cite[Theorem 1 and Theorem 2.6]{Mueller 1993}
who generalizes Cheeger's and his proof \cite{Cheeger (1979)}
and \cite{Mueller (1978)} of the Ray-Singer Conjecture
\mbox{$\rho_{\an}(M;V) = \rho_{\topological}(M;V)$}
for orthogonal representations $V$ to the unimodular
setting. Bismut and Zhang
 \cite[Theorem 0.1 and Theorem 0.2]{Bismut-Zhang (1992)}
have generalized M\"uller's work to all dimensions and
to a setting where $V$ is not necessarily unimodular.
We say more about their version of Theorem
\ref{equality of analytic and topological torsion}
in Section \ref{The fiber bundle formula for analytic torsion}.\par

Theorem \ref{equality of analytic and topological torsion}
and Theorem \ref{variation of Riemanninan metrics}
enables us to show that
Theorem \ref{fiber bundle formula for analytic torsion}
follows from its topological version. This will be done in
Section \ref{The fiber bundle formula for analytic torsion}
where we also explain how
Corollaries \ref{acyclic case},
\ref{spheres as fibers} and
\ref{spheres as base} follow from
Theorem \ref{fiber bundle formula for analytic torsion}
and its topological version.
In the topological case we can treat a more general setting.
Namely, we consider a fibration \mbox{$f : E \longrightarrow B$}
such that $B$ is a
connected finite $CW$-complex, the homotopy fiber has the
homotopy type of a finite $CW$-complex and a certain
cohomology class \mbox{$\theta_f \in H^1(B;\Wh(E))$} vanishes.
Then we get  for  a local coefficient system $V$ with unimodular
Hilbert structure, and fixed Hilbert structures on the
relevant singular  cohomology groups:\\[5mm]
{ \bf Theorem \ref{fibration formula}}
$$\rho(E;V) ~ = ~ \chi(B) \cdot \rho(F;V) +
\rho(B;\sum_{q} (-1)^q \cdot H^q_{\sing}(F;V)) + \rho^{\Serre}_{\sing}(f)
\qed$$

In particular the topological version includes manifolds
with boundary and the question whether Theorem
\ref{fiber bundle formula for analytic torsion} generalizes
to manifolds with boundary comes down to the question whether
Theorem \ref{equality of analytic and topological torsion}
and Theorem \ref{variation of Riemanninan metrics} generalize
to manifolds with boundary. At least under the condition that
$V$ is an orthogonal representation Cheeger and M\"uller's
results for closed manifolds have been extended to manifolds
with boundary (see \cite{Hassell (1994)},
\cite{Lott-Rothenberg (1991)} and
\cite{Lueck (1993)}). \par

There are also  interesting
generalizations of the topological versions to other settings.
For example one can consider an $L^2$-version. Or one
can substitute $V$ by a local
coefficient system of finitely generated projective modules
over a ring $R$. Then the torsion takes value in the
algebraic $K$-theory of $R$.

Instead of explicitely choosing Hilbert structures to define
torsion as a real number, one can consider it as an element in certain
determinant spaces. The latter approach has the advantage not to 
depend on artificial choices. However, it does not allow the generalizations
we have mentioned.
In particular, one is bound to the finite dimensional setting
(i.e.~even the singular cochain complex of a finite CW-complex is not allowed).
Also, if one carries out explicit computations,
in most cases it is preferable to deal with real numbers.
 Therefore we used the first method. In our context,
both languages are completely equivalent and we supply a dictionary to translate
between then in appendix \ref{Torsion and determinants}. One should also
mention that there is a third equivalent approach which uses norms constructed
on determinant lines (compare \cite{Bismut-Zhang (1992)}).

Finally we mention that in the case where $V$ is orthogonal
Dai and Melrose \cite{Dai-Melrose (1995)} have proven
the formula of Theorem
\ref{fiber bundle formula for analytic torsion}
in the adiabatic limit by completely different methods,
namely, by a careful analysis of the heat kernel in
the adiabatic limit and the adiabatic version of
the Leray-Serre spectral sequence \cite{Mazzeo-Melrose (1990)}.
We remark that in the special case
where $V$ is an orthogonal representation
\mbox{$\rho_{\an}(E;V)$} vanishes if \mbox{$\dim(E)$} is
even but this
is not true in general for a non-flat Riemannian metric on $V$. \par

The first two authors decidate this paper to their friend and colleague
Thomas Thielmann who died in a car accident in November 1994.

The paper is organized as follows\\[3mm]
\begin{tabular}{ll}
\ref{Simple structures on spaces}.
& Simple structures on spaces\\
\ref{Milnor torsion for cochain complexes}.
& Milnor torsion for cochain complexes\\
\ref{Milnor torsion for spaces and local coefficient systems}.
& Milnor torsion for spaces and local coefficient systems\\
\ref{Torsion and spectral sequences}.
& Torsion and spectral sequences\\
\ref{The fibration formula for Milnor torsion}.
& The fibration formula for Milnor torsion\\
\ref{The fiber bundle formula for analytic torsion}.
& The fiber bundle formula for analytic torsion\\
\ref{Torsion and determinants} & Torsion and determinants\\
 & References
\end{tabular}



\setcounter{section}{0}
\refstepcounter{section}
\typeout{------------------section 1 --------------------}

\tit{Simple structures on spaces}
\label{Simple structures on spaces}

In this section we explain additional structures on arbitrary topological
spaces and fibrations which allow  the definition of  torsion invariants on them. 
This extension
from the category of finite CW-complexes serves two purposes: on the one hand
there are interesting spaces which are not finite CW-complexes, f.i.\
classifiying spaces of discrete groups. On the other hand, our approach singles
out what exactly is used in the definition of torsion invariants, and this
definitely clarifies the exposition.

We will always assume for a pair of spaces $(X,A)$ that the inclusion
of $A$ into $X$ is a cofibration. This condition is satisfied if
$(X,A)$ is a pair of $CW$-complexes or if $X$ is a manifold with
submanifold $A$. A map
\mbox{$(F,f) : (X,A) \longrightarrow (Y,B)$}
of pairs is a {\em relative homotopy equivalence } if
\mbox{$F \cup \id: X \cup_f B \longrightarrow  Y$}
is a homotopy equivalence. Here
$X \cup_f B$  is obtained from $X \coprod B$
by identifying $a \in A$ with $f(a) \in B$. Notice
that $(F,f)$ is a homotopy equivalence of pairs if and only if $f$ is
a homotopy equivalence and $(F,f)$ is a relative homotopy equivalence.
Given a cellular relative  homotopy equivalence
\mbox{$(F,f) : (X,A) \longrightarrow (Y,B)$}
of pairs of finite $CW$-complexes, define its Whitehead torsion as the
Whitehead torsion of the homotopy equivalence
\mbox{$F \cup \id: X \cup_f B \longrightarrow  Y$} of finite $CW$-complexes
\begin{eqnarray}
\tau(F,f) & := &
\tau(F \cup \id) \hspace{10mm}
\in \Wh(Y).
\label{Whitehead torsion for homotopy equivalence of finite CW-complexes}
\end{eqnarray}
We refer to \cite[\S 6, \S 21 and \S 22]{Cohen (1973)}
for the definition of the
geometric Whitehead group and Whitehead
torsion and their identifications
with the algebraic Whitehead group $\Wh(\pi_1(Y))$ and Whitehead torsion.
Notice that $\tau(F,f)$ depends only
on the homotopy class of $(F,f)$
and satisfies the composition formula,
the formula for pairs and the product
formula as stated in
\ref{formulas for geometric torsion}.
This follows from  the special case $A  = \emptyset$
in \cite[\S 22 and \S 23]{Cohen (1973)}.
Given a pair $(Y,B)$, we call two relative homotopy equivalences
\mbox{$(F_i,f_i) : (X_i,A_i) \longrightarrow (Y,B)$} with pairs of finite
$CW$-complexes as source for $i =0,1$
equivalent if
\mbox{$\tau((G,g) \circ (F_0,f_0)) ~ = ~ \tau((G,g) \circ (F_1,f_1))$}
holds for any relative homotopy equivalence
\mbox{$(G,g): (Y,B)  \longrightarrow (Y^{\prime},B^{\prime})$}
into a pair of finite $CW$-complexes.

\begin{definition} \label{relative simple structure}
A {\em relative simple structure}  on a pair $(Y,B)$
is an equivalence class
of relative homotopy equivalences
\mbox{$(F,f) : (X,A) \longrightarrow (Y,B)$}
for a pair of finite $CW$-complexes as source.
\qed \end{definition}

\brn \label{representatives by homotopy equivalences of pairs}
Let $(X,A)$ be a pair with a relative simple structure and
\mbox{$g : A^{\prime\prime} \longrightarrow A$} be a homotopy equivalence
with a finite $CW$-complex as source. Then we can extend $g$ to
a representative for the simple structure
\mbox{$(G,g) :  (X^{\prime\prime},A^{\prime\prime}) \longrightarrow (X,A)$}
which is a homotopy equivalence of pairs as follows. \\
Choose some
representative
\mbox{$(F,f) : (X^{\prime},A^{\prime}) \longrightarrow (X,A)$}
for the given simple structure. Furthermore choose a homotopy
inverse \mbox{$g^{-1} : A \longrightarrow A^{\prime\prime}$}, a homotopy
$\phi$ from $g^{-1} \circ f$ to a cellular map
\mbox{$A^{\prime} \longrightarrow A^{\prime\prime}$}
and a homotopy $\psi$ from $g \circ g^{-1} \circ f$ to $f$.
Let $X^{\prime\prime}$ be the finite $CW$-complex
\mbox{$X^{\prime} \cup _{\phi_1} A^{\prime\prime}$}.
Define  \mbox{$(G,g) :  (X^{\prime\prime},A^{\prime\prime})
\longrightarrow (X,A)$}
by the following composition of (relative)
homotopy equivalences or their homotopy inverses
$$X^{\prime} \cup _{\phi_1} A^{\prime\prime} \longrightarrow
X^{\prime} \times [0,1] \cup _{\phi} A^{\prime\prime} \longleftarrow
X^{\prime} \cup _{g^{-1} \circ f} A^{\prime\prime}
\longrightarrow
X^{\prime} \cup_{g \circ g^{-1} \circ f} A$$
$$\longrightarrow  X^{\prime} \times [0,1] \cup_{\psi} A
\longleftarrow X^{\prime} \cup_{f}
A \longrightarrow X.$$
The proof that $(G,g)$ represents the given simple structure is done
by the results of \cite[\S 5]{Cohen (1973)}. \qed
\ern

\brn \label{simple structure on X from  A and (X,A)}
Given a (relative)  simple structure on $(X,A)$ and on $A$,
we construct a preferred simple structure on $X$ as follows.
Because of \ref{representatives by homotopy equivalences of pairs}
there is a homotopy equivalence of pairs
\mbox{$(G,g) :  (X^{\prime\prime},A^{\prime\prime}) \longrightarrow (X,A)$}
such that it represents the given relative simple structure on $(X,A)$
and $g$ represents the given simple structure on $A$. The preferred
simple structure on $X$ is then represented by $G$.
This is independent of the choice of $(G,g)$ by
\ref{formulas for geometric torsion}.
\qed\ern

Given a relative homotopy equivalence
\mbox{$(F,f) : (X,A) \longrightarrow (Y,B)$}
of pairs with relative simple
 structures,
we still can define its (relative) Whitehead torsion
\begin{eqnarray}
\tau(F,f) & \in & \Wh(Y)
\label{definition of geometric Whitehead torsion}
\end{eqnarray}
as follows. Choose representatives
\mbox{$(G,g) : (X^{\prime},A^{\prime}) \longrightarrow (X,A)$} and
\mbox{$(H,h) : (Y^{\prime},B^{\prime}) \longrightarrow (Y,B)$}
for the relative
structures. Because of
\ref{representatives by homotopy equivalences of pairs}  one can
arrange that $(H,h)$ is a homotopy equivalence of pairs. Define
$\tau(F,f)$ as the image of Whitehead torsion
\mbox{$\tau((H,h)^{-1} \circ (F,f) \circ (G,g))$} defined in
\ref{Whitehead torsion for homotopy equivalence of finite CW-complexes}
under the map
\mbox{$H_{\ast} : \Wh(Y^{\prime}) \longrightarrow \Wh(Y)$} induced by $H$.

\brn \label{formulas for geometric torsion}
We have already mentioned homotopy invariance, the composition formula
$$\tau((G,g) \circ (F,f)) = \tau(G,g) + G_{\ast}\tau(F,f)$$
the formula for pairs
$$\tau(F) = \tau(F,f) + i_{\ast}\tau(f)$$
and the product formula
$$\tau((F,f) \times \id_Y) = \chi(Y) \cdot i_{\ast}\tau(F,f) $$
where $i$ denotes the obvious inclusions. One easily checks
that they remain true in the more general case that
$(X,A)$ is not necessarily a pair of finite $CW$-complexes, but
carries a relative simple structure. \qed \ern

\brn \label{simple structure on E}
Let \mbox{$f : E \longrightarrow B$} be a fibration such that $B$ is a finite
$CW$-complex and the fiber has the
homotopy type of  a finite $CW$-complex.
Suppose that we are given a {\em cellular base point system}
\mbox{$\{b_c \mid c \in I\}$} for $B$, i.e.
a choice of points $b_c$ in the
interior $c^{\circ}$ for each $c \in I$
where here and elsewhere $I_n$ is
the set of $n$-cells and $I$ is the disjoint union of the $I_n$-s.
Furthermore suppose that we have specified
 a simple structure on each fiber
\mbox{$F_{b_c} = f^{-1}(b_c)$}. We want to
 define a simple structure
on $E$ depending only on these choices as follows.
Recall that any homotopy class
of paths $w$ from $b_0$ to $b_1$ defines a homotopy
class of homotopy equivalences
\mbox{$t_w: F_{b_0} \longrightarrow F_{b_1}$}
by the {\em fiber transport} \cite[15.12]{Switzer (1975)}.\par

Let $E_n$ be $f^{-1}(B_n)$. As \mbox{$B_{n-1} \longrightarrow B_n$} is a
cofibration,
the same is true for \mbox{$E_{n-1} \longrightarrow E_n$}
\cite[I.7.14]{Whitehead (1978)}. Because of the construction
\ref{simple structure on X from  A and (X,A)} it suffices to specify
a relative simple structure for $(E_n,E_{n-1})$ for all $n \ge 0$.
 This will be
done by  the next Lemma
\ref{simple structure on (E_n,E_(n-1)} taking into account
\ref{formulas for geometric torsion} and
that the Whitehead torsion of any
homotopy equivalence \mbox{$(D^n,S^{n-1})
\longrightarrow (D^n,S^{n-1})$} is trivial
since $D^n$ is simply-connected and
hence $\Wh(D^n)$ is trivial. \qed \ern

\begin{lemma} \label{simple structure on (E_n,E_(n-1)}
Suppose we have specified a cellular base point system for $B$ and
for each element in $I_n$ an orientation.
Then there is a relative homotopy equivalence which is uniquely defined
up to homotopy
$$\coprod_{c \in I_n} F_{b_c} \times (D^n,S^{n-1})
\longrightarrow (E_n,E_{n-1}).$$
If we change the orientation of  the cell $c$, then the map is changed
by the selfhomotopy equivalence
\mbox{$\id \times s : F_{b_c} \times (D^n,S^{n-1})
\longrightarrow F_{b_c} \times (D^n,S^{n-1})$} where $s$ is a map of
degree $-1$. If we change the base point $b_c$ of $c$ to $b^{\prime}_c$,
then the map is changed by the homotopy equivalence
\mbox{$t \times \id : F_{b_c} \times (D^n,S^{n-1})
\longrightarrow
F_{b_c^{\prime}} \times (D^n,S^{n-1})$} where
\mbox{$t : F_{b_c} \longrightarrow F_{b^{\prime}_c}$} is given
by the fiber transport
along any path in $c^\circ$ connecting $b_c$ and $b^{\prime}_c$.
\end{lemma}
\proof The map will be constructed as the
composition of the following four
 maps or their homotopy inverses.\par

There is up to homotopy one orientation preserving
homotopy equivalence of pairs
\mbox{$(D^n,S^{n-1}) \longrightarrow (c^{\circ},c^{\circ}-b_c)$}.
This gives the first map
$$\coprod_{c \in I_n} F_{b_c} \times (D^n,S^{n-1}) \longrightarrow
\coprod_{c \in I_n} F_{b_c} \times (c^{\circ},c^{\circ}-b_c).$$
Choose a homotopy
\mbox{$\phi_c : c^{\circ} \times [0,1] \longrightarrow B$}
from the canonical inclusion $c^{\circ} \longrightarrow B$
to the constant map
with value $b_c$ such that its evaluation at $b_c$ gives a path within
$c^{\circ}$. By the homotopy lifting property we obtain a strong
fiber homotopy equivalence which is unique up to fiber homotopy
\cite[Proposition 15.11]{Switzer (1975)}
$$F_{b_c} \times (c^{\circ},c^{\circ} - b_c) \longrightarrow
(E|_{c^{\circ}},E|_{c^{\circ}-b_c}).$$
The second map is the disjoint union of these maps over $I_n$.
The third map is the relative homotopy
equivalence given by the inclusion
$$\coprod_{c \in I_n} (E|_{c^{\circ}},E|_{c^{\circ}-b_c})
\longrightarrow
(E_n,E|_{B_n - \{ \coprod_{c \in I_n}  b_c\}}).$$
The fourth map is the  homotopy equivalence of pairs given by the inclusion
$$(E_n,E_{n-1}) \longrightarrow
(E_n,E|_{B_n - \{ \coprod_{c \in I_n}  b_c\}}).$$
This finishes the proof of Lemma \ref{simple structure on (E_n,E_(n-1)}.
\qed

\brn \label{change of simple structure under change of base point system}
Let \mbox{$\{b_c \mid c \in I\}$} and
\mbox{$\{b_c^{\prime} \mid c \in I\}$}
be two base point systems and suppose
$F_{b_c}$ and $F_{b^{\prime}_c}$ come
with  simple structures. Let $\sigma$ and $\sigma^{\prime}$
be the simple structures
on $E$ given by the construction
\ref{simple structure on E} for these two choices. Let
\mbox{$t_c : F_{b_c} \longrightarrow F_{b^{\prime}_c}$}
be the homotopy equivalence
which is given by the fiber transport along any
path in $c^{\circ}$ from
$b_c$ to $b_c^{\prime}$. Denote by
\mbox{$i(b^{\prime}_c) : F_{b^{\prime}_c} \longrightarrow E$} the inclusion.
Then we get from  \ref{formulas for geometric torsion} and
Lemma \ref{simple structure on (E_n,E_(n-1)}
$$\tau(\id : (E,\sigma) \longrightarrow (E,\sigma^{\prime})) ~ = ~
\sum_{n \ge 0} (-1)^n \cdot
\sum_{c \in I_n} i(b^{\prime}_c)_{\ast}\tau(t_c : F_{b_c}
\longrightarrow F_{b_c^{\prime}}). \qed$$
\ern

Next we give a criterion when the choice
of base point system does not affect the simple structure on $E$.
Fix a base point $b \in B$. Given an element $w$ in $\pi_1(B,b)$,
define
$$\theta_f(w) ~ = ~ i(b)_{\ast}\tau(t_w : F_b \longrightarrow F_b) ~
\in \Wh(E)$$
where $t_w$ is the fiber transport and
\mbox{$i(b) : F_b \longrightarrow E$} is the inclusion. As
\mbox{$i(w(1)) \circ t_w \simeq i(w(0))$}
holds for any path $w$ in $B$, one easily checks using
\ref{formulas for geometric torsion} that this defines a homomorphism
from $\pi_1(B,b)$ to $\Wh(E)$ and thus a cohomology class which is
independent of the choice of $b \in B$
\begin{eqnarray}
\theta_f & \in & H^1(B;\Wh(E)).
\label{cohomology class of fibration}
\end{eqnarray}

\begin{definition} \label{coherent choice of simple structures on the fibers}
A {\em choice of simple structures on the fibers}
 is a choice of simple structure on each fiber $F_b$. It is called
{\em coherent} if for any path $w$ in $B$ we get
$$i(w(1))_{\ast}\tau(t_w: F_{w(0)} \longrightarrow
F_{w(1)}) ~ = ~ 0 ~ \in \Wh(E). \qed$$
\end{definition}

\brn \label{simple structure on E for coherent choices}
Notice that a coherent choice of simple structures on the fibers exists
if and only if $\theta_f$ is trivial.
If we fix a coherent choice of simple
structures
\mbox{$\{\sigma(F_b) \mid b \in B\}$} on the fibers,
the induced simple structure $\sigma$ on $E$
of \ref{simple structure on E}
is independent of the cellular base point system by
\ref{change of simple structure under change of base point system}.
Moreover, if we have a different choice of
coherent structures on the fibers, then
$$\tau(\id : (E,\sigma) \longrightarrow (E,\sigma^{\prime})) ~ = ~
\chi(B) \cdot
i(b)_{\ast}(\tau(\id : (F_b,\sigma(F_b)) \longrightarrow
(F_b,\sigma^{\prime}(F_b))).$$
In particular we see for a fibration
$p : E \longrightarrow B$ over a finite
$CW$-complex $B$ such that the homotopy fiber has the homotopy
type of a finite $CW$-complex  and $\chi(B) =0$
 and $\theta_f =0$ holds that
$E$ has a preferred simple structure.
\qed \ern

\begin{lemma} \label{smooth bundles and simple structures}
Let \mbox{$F \longrightarrow E \xrightarrow{f} B$}
be a smooth bundle of compact smooth
manifolds. Equip $B$ and each fiber with the simple structure given
by a  smooth triangulation. This gives a coherent choice of simple
structures on the fibers. Then the simple structure on $E$
given by a smooth triangulation agrees with the one given by
\ref{simple structure on E}. \qed \end{lemma}

\begin{remark} \label{simple structure on B suffices} \em
The construction \ref{simple structure on E}
can be extended
to the case where $B$ is not necessarily a $CW$-complex but
carries a simple structure. Namely, choose a representative
\mbox{$g : X \longrightarrow B$} of the simple structure.
The pull back construction yields a fibration
\mbox{$\overline{f} : g^{\ast}E \longrightarrow X$}
and a fiber homotopy equivalence
\mbox{$\overline{g} : g^{\ast}E \longrightarrow E$}. Equip $\overline{f}$
with the coherent
choice of simple structures on the fibers induced by $\overline{g}$
and the given one on $f$. Construction
\ref{simple structure on E for coherent choices}
applies to $\overline{f}$
and gives a simple structure on $g^{\ast}E$.  Equip $E$ with
the simple structure for which $\tau(\overline{g})$ vanishes.
It is not hard to check that this is independent of the choice
of the representative $g$. The main step is to show in the case
where $X$ and $B$ are finite $CW$-complexes and $g$ is an
elementary expansion that $\tau(\overline{g})$ vanishes with
respect to the simple structures on $E$ and $g^{\ast}E$ given
by construction \ref{simple structure on E for coherent choices}.
\qed\em
\end{remark}



\refstepcounter{section}
\typeout{----------  section 2   ------------------}

\tit{Milnor torsion for cochain complexes}
\label{Milnor torsion for cochain complexes}

In this section we give a brief introduction to torsion invariants of
cochain complexes as defined in the introduction.  We give no proofs but refer to \cite{Milnor (1966)}
and \cite {Lueck-Rothenberg (1991)}.\par
Given linear isomorphism \mbox{$f$}, $g$ of
finite-dimensional real Hilbert spaces, the number $[[\cdot]]$ of \ref{[[]]} has the following properties:
\begin{eqnarray}
[[f]] & = & \ln(|\det(A)|) ~ \in \rr
\end{eqnarray}
where $A$ is the matrix describing $f$ with respect
to some choice of orthonormal basis for source and range.
\begin{eqnarray*}
[[f \circ g]] & = & [[f]] + [[g]];
\\
\left[\left[\squarematrix{f}{h}{0}{g}\right]\right] & = &
[[f]] + [[g]] ;
\\
\mbox{$[[f^{\ast}]]$} & = & [[f]].
\end{eqnarray*}
Let $C = C^{\ast}$ be a finite Hilbert cochain complex, i.e. a cochain
complex of finite-dimensional Hilbert spaces such that
\mbox{$C^i = 0$} for \mbox{$|i| \ge N$} for some natural number $N$. If
$C$ is acyclic, we defined in \ref{Reidemeistert torsion for acyclic chain complexes}
$ \rho(C)  := 
[[(c^{\ast} + \gamma^{\ast}) : C^{\ev} \longrightarrow
C^{\odd}]] ~ \in \rr $,
where $c^{\ast}$ is the differential and $\gamma^{\ast}$ a
chain contraction. This is independent of the choice of $\gamma^{\ast}$.
For \mbox{$f : C \longrightarrow D$} a chain homotopy
equivalence of finite Hilbert cochain complexes, $t(f)$ was defined in
\ref{torsion of chain homotopy equivalence}.
It turns out that $t(f)$ depends only on the homotopy class of
$f$ and
\mbox{$t(f \circ g) = t(f) + t(g)$}. Notice that with these conventions
we get for an isomorphisms \mbox{$f : C \longrightarrow D$}
of finite Hilbert cochain complexes
$$t(f) ~ = ~ \sum_{n} (-1)^n \cdot [[f^n]].$$
Let $C$ be a finite Hilbert cochain complex such that
$H(C^{\ast})$ carries a Hilbert structure, i.e. $H^n(C^{\ast})$
is equipped with a Hilbert space structure for each $n \in \zz$.
If
\mbox{$i : H(C) \longrightarrow C$} is the chain map which induces on
cohomology the indentity, we defined
$\rho(C) := -t(i) ~ \in \rr$.\\
The minus sign ensures
that this definition coincides with the one in
\ref{Reidemeistert torsion for acyclic chain complexes}
in the acyclic case.
If we fix an orthonormal basis for each $C^n$ and each $H^n(C)$,
then the logarithm of the torsion defined by Milnor
\cite[page 365]{Milnor (1966)} is $\rho(C)$. \par

Let $C$ be a finite Hilbert cochain complex and
equip $\ker(c^p)$ and $\im(c^{p-1})$ with the Hilbert substructures
and $H^p(C) = \ker(c^p)/\im(c^{p-1})$ and $C^p/\ker(c^p)$
with the quotient Hilbert structures. If
\mbox{$\overline{c^p} : C^p/\ker(c^p)\longrightarrow \im(c^p)$}
is the obvious isomorphism induced by $c^p$ and
\mbox{$\overline{\Delta^{p}} : C^p/\ker(\Delta^p) \longrightarrow
C^p/\ker(\Delta^p)$}
is the automorphism induced by the endomorphism
\mbox{$\Delta^p = (c^p)^{\ast}\circ c^p + c^{p-1}\circ (c^{p-1})^{\ast}:
C^p \longrightarrow C^p$}, then we get
\begin{eqnarray}
\rho(C) & = & \sum_{p} (-1)^p \cdot [[\overline{c^p}]]
~ =~  - \sum_{p} (-1)^p \cdot p \cdot
\frac{1}{2} \cdot \ln\left(\det(\overline{\Delta^p})\right).
\label{special Hilbert structure on cohomology}
\end{eqnarray}

A {\em simple structure} on a (real) cochain complex $C$ is an
equivalence class of chain homotopy equivalences
\mbox{$u : \overline{C} \longrightarrow C$}
with a finite Hilbert cochain complex as source where $u$ and
\mbox{$v : \overline{\overline{C}} \longrightarrow C$}
are  equivalent if
\mbox{$t(v^{-1} \circ u)$} vanishes. Let
\mbox{$f : C \longrightarrow D$} be a chain homotopy equivalence
of cochain complexes with simple structure. Define
\begin{eqnarray}
t(f) & = & t(v^{-1} \circ f \circ u) ~  \in \rr
\label{torsion for simple chain structures}
\end{eqnarray}
for any representatives
\mbox{$u : \overline{C} \longrightarrow C$}
and \mbox{$v : \overline{D} \longrightarrow D$}
of the simple structures. Let $C$ be a
cochain complex with simple structure
such that $H(C)$ has a Hilbert structure. Define
\begin{eqnarray}
\rho(C) & := & \rho(\overline{C}) ~ \in \rr
\label{Milnor torsion for simple chain structures}
\end{eqnarray}
for any representative \mbox{$u : \overline{C} \longrightarrow C$}
where we use the Hilbert structure on $H(\overline{C})$ for
which $H(u)$ is an isometry.\par

Next we collect the basic properties of these invariants. We mention
that one firstly verifies them for finite Hilbert cochain complexes
and uses this to show that the definitions and results extend to
cochain complexes with simple structures. \par

If \mbox{$\shortexactsequence{C}{}{D}{}{E}$} is an exact sequence
of finite Hilbert cochain complexes, we can view
\mbox{$C^n \longrightarrow D^n \longrightarrow E^n$}
as an acyclic finite Hilbert cochain complex concentrated
in dimension $0$, $1$ and $2$ and define
$$\rho(C \longrightarrow D \longrightarrow E) ~:= ~
\sum_{n \in \zz} (-1)^n \cdot \rho(C^n
\longrightarrow D^n \longrightarrow E^n).$$
This extends to an exact sequence
\mbox{$\shortexactsequence{C}{}{D}{}{E}$}
of cochain complexes with simple structures as follows. Construct a
commutative diagram of cochain complexes
\begin{center}
$\begin{CD}
0 @>>>  C  @>>>  D  @>>>  E   @>>> 0\\
 & & @Af AA  @Ag AA  @Ah AA  \\
0  @>>>  \overline{C}    @>>>  \overline{D}     @>>>
 \overline{E}    @>>>&0
\end{CD}$
\end{center}
with the property that the rows are exact and the vertical
arrows are homotopy equivalences which represent the given
simple structures. In particular the lower row is an exact
sequence of finite Hilbert cochain complexes and we put
\begin{eqnarray}
\rho(C \longrightarrow D \longrightarrow E) & := &
\rho(\overline{C} \longrightarrow \overline{D}
\longrightarrow \overline{E}).
\label{torsion of exact sequence}
\end{eqnarray}

\brn \label{main properties of torsion for chain complexes}
In the following list of basic properties
 all cochain complexes come with simple structures.

\begin{enumerate}
\item homotopy invariance\\
$f \simeq g ~ \Rightarrow ~ t(f) = t(g)$;

\item composition formula
$$t(f \circ g) = t(f) + t(g);$$

\item exactness\\
Given a commutative diagram of the shape above
with exact rows and homotopy equivalences as vertical arrows, then
$$t(f) - t(g) + t(h) =
\rho(\overline{C} \longrightarrow \overline{D} \longrightarrow \overline{E}) -
\rho(C \longrightarrow D \longrightarrow E);$$

\item sum formula\\
Let \mbox{$\shortexactsequence{C}{}{D}{}{E}$} be exact and
the cohomology of $C$, $D$ and $E$ come with Hilbert structures.
Let $L\!H\!S$ be the acyclic finite Hilbert cochain complex given by
the long  cohomology sequence where
\mbox{$L\!H\!S^0 = H^0(C)$}.
Then:
$$\rho(D) - \rho(C) - \rho(E) ~ = ~
\rho(L\!H\!S) - \rho(C \longrightarrow D \longrightarrow E);$$

\item transformation formula\\
If \mbox{$f : C \longrightarrow D$} is a
homotopy equivalence and the cohomology
of $C$ and $D$ come with Hilbert structures, then:
$$\rho(C)- \rho(D) ~ = ~ t(f) -
\sum_{n \in \zz} (-1)^n \cdot [[H^n(f)]].$$
\end{enumerate}
\ern



\refstepcounter{section}
\typeout{------------------ section 3 ---------------------}

\tit{Milnor torsion for spaces and local coefficient systems}
\label{Milnor torsion for spaces and local coefficient systems}
The fundamental groupoid $\Pi(X)$
of a space $X$ has as objects points in $X$ and
a morphism \mbox{$w : x \longrightarrow y$}
is a homotopy class relative end
points of paths in $X$ from $y$ to $x$. For $x \in X$ let
$\widetilde{X}(x)$ be the set of all morphisms
\mbox{$w : x \longrightarrow y$} with $x$ as source.
The projection
\mbox{$p(x) : \widetilde{X}(x) \longrightarrow X$} sends $w$ to $y$.
Assume that
$X$ is locally path-connected and semi-locally simply connected.
This condition is always satisfied if $X$ is a $CW$-complex or manifold.
It ensures that there is precisely one topology on $\widetilde{X}(x)$
such that $p(x)$ is a model for the universal covering of the path
component of $X$ containing $x$.
Thus we get a contravariant functor
$$\widetilde{X} : \Pi(X) \longrightarrow \SPACES.$$
Composing it with the covariant functor singular chain complex with
real  coefficients yields the contravariant functor
$$C_{\ast}^{\sing}(\widetilde{X}) : \Pi(X)
\longrightarrow \rr-\CHAIN$$
A local coefficient system $V$ on $X$ is a contravariant functor
\mbox{$V : \Pi(X) \longrightarrow \rr-\VECTOR$} into the category
of finite-dimensional real vector spaces. For instance a flat
vector bundle over $X$ defines a local coefficient system.
Define the
{\em singular cochain complex and the singular cohomology of $X$
with coefficient in $V$} by
\begin{eqnarray}
C^{\ast}_{\sing}(X;V) & := &
\hom(C^{\sing}_{\ast}(\widetilde{X}),V)
\label{cochain complex and cohomology with coefficients}
\\
H^{\ast}_{\sing}(X;V) & := &
H^{\ast}(C^{\ast}_{\sing}(X;V))
\nonumber
\end{eqnarray}
where $\hom$ denotes the real vector space
of natural transformations. \par

Suppose that $X$ is a $CW$-complex. Then $\widetilde{X}$ becomes a
functor from $\Pi(X)$ into the category of $CW$-complexes and
we can define the cellular versions
\mbox{$C^{\ast}_{\cell}(X;V)$} and \mbox{$H^{\ast}_{\cell}(X;V)$}
of the definitions
above by using the cellular chain complex instead of the singular
one. There is a homotopy equivalence
\begin{eqnarray}
C^{\ast}_{\cell}(X;V)
& \longrightarrow &
C^{\ast}_{\sing}(X;V)
\label{comparing singular and cellular}
\end{eqnarray}
which is unique up to homotopy and natural in $X$ and $V$
\cite[page 263]{Lueck (1989)}. It induces a natural isomorphism
\mbox{$H^{\ast}_{\cell}(X;V)
\longrightarrow H^{\ast}_{\sing}(X;V)$}.\par

Let $X$ be a space and $V^0$, $V^1$, $\ldots $, $V^r$ local
coefficient systems. We say
\mbox{$\sum_{q=0}^r (-1)^q \cdot V^q$} is
{\em unimodular} if  we have for each automorphism
\mbox{$w : x \longrightarrow x$} in $\Pi(X)$
$$\prod_{q=0}^r \left|\det(V_w^q : V_x^q
\longrightarrow V_x^q)\right| ^{(-1)^q}
  ~  = ~ 1.$$
A {\em Hilbert structure} on \mbox{$\sum_{q=0}^r (-1)^q \cdot V^q$}
is a choice of Hilbert space structure on $V^q_x$ for each
\mbox{$q \in \{1,2 \ldots r\}$} and $x \in X$. It is called
{\em unimodular} if we have
for each morphism \mbox{$w : y \longrightarrow x$} in $\Pi(X)$
$$\sum_{q = 0}^r (-1)^q \cdot [[V_w^q : V_x^q
\longrightarrow V_y^q]] ~=~ 0$$
in the notation of \ref{[[]]}. Notice that
\mbox{$\sum_{q=0}^r (-1)^q \cdot V^q$} is unimodular
if and only if admits a unimodular Hilbert structure. \par

\begin{definition} \label{Milnor torsion}
Let $X$ be a space and $V^0$, $V^1$, $\ldots $, $V^r$ be local
coefficient systems. Assume that $X$ comes with a simple structure,
\mbox{$\sum_{q=0}^r (-1)^q \cdot V^q$} with a unimodular Hilbert
structure and $H^{\ast}_{\sing}(X;V)$ with a Hilbert structure. Define
the {\em Milnor torsion}
$$\rho(X;\sum_{q=0}^r (-1)^q \cdot V^q) ~ := ~
\sum_{q =0}^r (-1)^q \cdot \rho(C_{\sing}^{\ast}(X;V^q))$$
where $\rho(C_{\sing}^{\ast}(X;V^q))$ was defined in
\ref{Milnor torsion for simple chain structures} and we use the
simple structure on $C_{\sing}^{\ast}(X;V^q)$ defined below. \qed
\end{definition}
Let \mbox{$f : Y \longrightarrow X$} be a representative for the simple
structure on $X$. Fix a cellular base point system
\mbox{$\{b_c \mid c \in I\}$}, i.e. a choice of base points
\mbox{$b_c \in c^{\circ}$}  for each cell $c$ of $Y$, and an
orientation for each cell $c$. For each cell $c$ there is
precisely one lift
\mbox{$\tilde{c} \subset \widetilde{Y}(b_c)$} which contains the
canonical base point in $\widetilde{Y}(b_c)$ given by the constant
path at $b_c$. The orientation of $c$ induces an orientation of
$\tilde{c}$. Thus we get an element
\mbox{$[\tilde{c}] \in C^{\cell}_{\dim(c)}(\widetilde{Y}(b_c))$}.
Define an isomorphism
\begin{eqnarray}
C^n_{\cell}(Y,f^{\ast}V^q)\longrightarrow
\oplus_{c \in I_n} V_{f(b_c)}
& \hspace{10mm} &
t \mapsto \left(t(b_c)([\tilde{c}])\right)_{c \in I_n}.
\label{identification of modules}
\end{eqnarray}
Equip $C^n_{\cell}(Y,f^{\ast}V^q)$ with the Hilbert structure induced
by the isomorphism \ref{identification of modules} above and the given
Hilbert space structure on the various $V_{f(b_c)}$-s. The desired
simple structure on $C_{\sing}^{\ast}(X;V^q)$ is represented by the
composition of the chain homotopy equivalences given by
\ref{comparing singular and cellular} and $f$
$$C_{\cell}^{\ast}(Y;f^{\ast}V^q)
\longrightarrow C_{\sing}^{\ast}(Y;f^{\ast}V^q)
\xrightarrow{(f^{\ast})^{-1}} C_{\sing}^{\ast}(X;V^q).$$

The choice of the orientations of the cells does not affect this simple
structure. If we change the base point system, the simple structure and
hence \mbox{$\rho(C_{\sing}^{\ast}(X;V^q))$} changes. However,
\mbox{$\sum_{q =0}^r (-1)^q \cdot \rho(C_{\sing}(X;V))$} does not change
because of the formulas
\ref{main properties of torsion for chain complexes}
and the assumption that the Hilbert structure on
\mbox{$\sum_{q=0}^r (-1)^q \cdot V^q$} is unimodular. \par

If \mbox{$\sum_{q=0}^r (-1)^q \cdot V^q$} is unimodular,
we can define a homomorphism
\begin{eqnarray}
\Phi(\sum_{q=0}^r (-1)^q \cdot V^q) : \Wh(X)
& \longrightarrow & \rr
\label{homomorphism from Wh to R}
\end{eqnarray}
as follows. An element $[u] \in \Wh (X)$ can
be represented by an automorphism
\mbox{$u : M \longrightarrow M$} of a finitely generated free
$\zz\Pi(X)$-module $M$.
Composition with $u$ defines an automorphism
\mbox{$u^q : \hom(M,V^p) \longrightarrow \hom(M,V^q)$}
of a finite-dimensional real
vector space. Put
$$\Phi(\sum_{q=0}^r (-1)^q \cdot V^q)([u]) := \sum_{q=0}^r (-1)^q \cdot
 \ln|\det(u^q)|.$$
The unimodularity condition on \mbox{$\sum_{q=0}^r (-1)^q \cdot V^q$}
ensures that trivial units are send to $0$. The following lemma is a
consequence of the definitions and the formulas
\ref{main properties of torsion for chain complexes}
\begin{lemma} \label{comparing Whitehead and Milnor torsion} Let
\mbox{$f : X \longrightarrow Y$} be a homotopy equivalence of spaces with
simple structures. Let  $V^0$, $V^1$, $\ldots $, $V^r$ be local
coefficient systems on $Y$. Assume that
\mbox{$\sum_{q=0}^r (-1)^q \cdot V^q$} comes with a unimodular
Hilbert
structure and $H^{\ast}_{\sing}(X;f^{\ast}V)$ and
$H^{\ast}_{\sing}(Y;V)$
come with Hilbert structures. Then:
\begin{eqnarray*}
\rho(Y;V) - \rho(X;f^{\ast}V) & = &
\Phi(\sum_{q=0}^r (-1)^q \cdot V^q)(\tau(f))
\\
& &  -
\sum_{q = 0}^r (-1)^q \cdot \sum_{p\ge 0} (-1)^p \cdot
[[f^p : H^p_{\sing}(Y;V^q) \longrightarrow
H^p_{\sing}(X;f^{\ast}V^q)]]. \hspace{3mm} \rule{2mm}{2mm}
\end{eqnarray*}
\end{lemma}


\refstepcounter{section}
\typeout{----------- section 4 -----------------------}

\tit{Torsion and spectral sequences}
\label{Torsion and spectral sequences}

Let $C$ be a cochain complex with finite descending filtration
by cochain complexes
$F_pC$
$$C= F_{0}C \supset F_1C \supset F_2C \supset \ldots \supset
F_pC \supset F_{p+1}C \supset \ldots.$$
Finite means that there is a natural number $l$ with
\mbox{$F_lC = \{0\}$}. Put \mbox{$F_pC = C$} for \mbox{$p \le -1$}.
We recall the construction of the
associated spectral cohomology sequence
\mbox{$(E_*^{*,*},d_*^{*,*})$} converging to $H^{\ast}(C)$
since we will need it explicitely
(see \cite{Cartan-Eilenberg (1956)} or
\cite{McCleary (1985)}).
In the sequel
$\partial$ denotes the boundary operator of the long exact
cohomology sequence associated to a short exact sequence of
cochain complexes.  We will abbreviate $F_pC$ by $F_p$.
Define  for $r \ge 1$:
\begin{eqnarray*}
Z_r^{p,q} & := &  \im(H^{p+q}(F_p/F_{p+r}) \longrightarrow
H^{p+q}(F_p/F_{p+1}))
\\
B_r^{p,q} & := & \im(H^{p+q-1}(F_{p-r+1}/F_p)
\xrightarrow{\partial} H^{p+q}(F_p/F_{p+1}))
\\
Z_{\infty}^{p,q} & := & \im(H^{p+q}(F_p) \longrightarrow
H^{p+q}(F_p/F_{p+1}))
\\
B_{\infty}^{p,q} & := & \im(H^{p+q-1}(C/F_p)
\xrightarrow{\partial}  H^{p+q}(F_p/F_{p+1}))
\\
E_r^{p,q} & := & Z_r^{p,q}/B_r^{p,q}
\\
E_{\infty}^{p,q} & := & Z_{\infty}^{p,q}/B_{\infty}^{p,q}
\\
F^{p,q} & := & \im(H^{p+q}(F_p) \longrightarrow H^{p+q}(C))
\end{eqnarray*}
We have the inclusions:
$$\{0\} = B_1^{p,q} \subset \ldots \subset B_r^{p,q} \subset
B_{\infty}^{p,q} \subset Z_{\infty}^{p,q} \subset Z_r^{p,q} \subset
\ldots \subset Z_1^{p,q} = H^{p+q}(F_p/F_{p+1})$$
The map
\mbox{$H^{p+q}(F_p/F_{p+r}) \longrightarrow
H^{p+q}(F_p/F_{p+1})$} resp.\
\mbox{$H^{p+q}(F_p/F_{p+r}) \longrightarrow
H^{p+q+1}(F_{p+r}/F_{p+r+1})$} is induced by the inclusion resp.\
is a boundary operator.
We get epimorphisms
$$H^{p+q}(F_p/F_{p+r}) \longrightarrow Z_r^{p,q}/Z_{r+1}^{p,q}$$
and
$$H^{p+q}(F_p/F_{p+r}) \longrightarrow
B_{r+1}^{p+r,q-r+1}/B_r^{p+r,q-r+1}.$$
Now the standard diagram chase shows that
these maps have the same kernel.
Hence we  obtain canonical isomorphisms
$$\gamma_r^{p,q} : Z_r^{p,q}/Z_{r+1}^{p,q}
\longrightarrow B_{r+1}^{p+r,q-r+1}/B_r^{p+r,q-r+1}
\hspace{10mm} \mbox{ for } r \ge 1.$$
Analogously one gets a natural isomorphism
$$\psi^{p,q} : F^{p,q}/F^{p+1,q-1} \longrightarrow E_{\infty}^{p,q}.$$
We define the differential
$$d_r^{p,q} : E_r^{p,q} = Z_r^{p,q}/B_r^{p,q} \longrightarrow
E_r^{p+r,q-r+1} = Z_r^{p+r,q-r+1}/B_r^{p+r,q-r+1}$$
by the composition:
$$Z_r^{p,q}/B_r^{p,q} \longrightarrow Z_r^{p,q}/Z_{r+1}^{p,q}
\xrightarrow{\gamma_r^{p,q}}
B_{r+1}^{p+r,q-r+1}/B_r^{p+r,q-r+1} \longrightarrow
Z_r^{p+r,q-r+1}/B_r^{p+r,q-r+1}$$
We obtain cochain complexes
$$(E_r^{p+r\ast,q-(r-1)\ast},d_r^{p+r\ast,q-(r-1)\ast})
\hspace{10mm} r \ge 0$$
if we use for $r \ge 1$ the definition above and
define \mbox{$E_0^{p,\ast}$} to be the $(-p)$-th suspension
\mbox{$\Sigma^{-p} F_p/F_{p+1}$} of \mbox{$F_p/F_{p+1}$}.
Since \mbox{$\ker(d_r^{p,q})$} is $Z_{r+1}^{p,q}/B_r^{p,q}$ and
\mbox{$\im(d_r^{p,q})$} is
\mbox{$B_{r+1}^{p+r,q-r+1}/B_r^{p+r,q-r+1}$},
we obtain a canonical
isomorphism
$$\phi_r^{p,q} : H^{0}(E_r^{p+r\ast,q-(r-1)\ast}) \longrightarrow
 E_{r+1}^{p,q} \hspace{10mm} \mbox{ for } r \ge 0$$
in the case $r \ge 1$, in  the case $r = 0$
we use the identification
\mbox{$H^0(E_0^{p,q+\ast}) = H^{p+q}(F_p/F_{p+1})$}.

\brn \label{data for chain complexes}
Now suppose we have fixed the following data
\begin{enumerate}
\item Simple structures on
\mbox{$F_p/F_{p+1}$};

\item a Hilbert structure on $H(C)$. \qed
\end{enumerate}
\ern

Equip $F_p$ inductively with the simple structure for which
the number defined in \ref{torsion of exact sequence} satisfies
$$\rho(F_{p+1} \longrightarrow F_p \longrightarrow F_p/F_{p+1}) ~=~ 0.$$
In particular we get a preferred simple
structure on $C= F_0$ and $\rho(C)$ is defined.
Equip $F^{p,q} \subset H^{p+q}(C)$ with the Hilbert substructure
and $F^{p,q}/F^{p+1,q-1}$ with the Hilbert quotient structure. Do iteratively 
the same for $E^{p,q}_r$ and $H(E^{p+r*,q-(r-1)*}_r)$.
Observe that $\phi^{p,q}_r$ become isometries.
Notice that \ref{data for chain complexes}.1 implies that
$E^{p,q}_r$ is finite-dimensional for $r \ge 1$.

\begin{definition} \label{torsion of a filtered chain complex}
Define for the finite descending filtration $F_{\ast}C$
with respect to the data \ref{data for chain complexes} and
choices above :
\begin{eqnarray*}
\rho_{\fil}^{\ge 2}(F_{\ast} C) & := &
\sum_{r \ge 2} \sum_{p=0}^{r-1} \sum_{q } (-1)^{p+q}
\cdot \rho(E_r^{p+r\ast,q-(r-1)\ast})
~ - ~ \sum_{p,q} (-1)^{p+q} \cdot [[\psi^{p,q}]];
\\[3mm]
\rho_{\fil}(F_{\ast}C) & := &
\sum_{p} \rho(F_p/F_{p+1}) ~ + ~
\sum_{r \ge 1} \sum_{p=0}^{r-1} \sum_{q } (-1)^{p+q}
\cdot \rho(E_r^{p+r\ast,q-(r-1)\ast})
\\
& &  - ~ \sum_{p,q} (-1)^{p+q} \cdot [[\psi^{p,q}]];
\end{eqnarray*}
where $\rho$ was defined in
\ref{Milnor torsion for simple chain structures}
and $[[\hspace{1mm}]]$ in \ref{[[]]}. \qed
\end{definition}

\begin{remark}\label{def of rho ge 2}
$\rho_{\fil}(F_{\ast}C)$ depends on the choices 
\ref{data for chain complexes}.\\
For the definition of $\rho_{\fil}^{\ge 2}(F_{\ast} C$ one can fix Hilbert
structures on $E^{p,q}_2$ instead of data \ref{data for chain complexes}.1 and
equip $E^{p,q}_r$ and $H(E^{p+r*,q-(r-1)*}_r$ ($r\ge 2$) with the corresponding Hilbert sub- and quotient structures. Then, $\rho_{\fil}^{\ge 2}(F_{\ast} C$ depends not on \ref{data for chain complexes}.1 but on these choices  as follows:
 if $U_2^{p,q}: \overline{E}_2^{p,q}\to E_2^{p,q}$ is the identity
on $E_2$ equipped with two different Hilbert structures, then
\[ \overline{\rho_{\fil}^{\ge 2}(F_{\ast} C)}-\rho_{\fil}^{\ge 2}(F_{\ast} C)=
\sum_{p,q}(-1)^{p+q}[[U_2^{p,q}]] \] 
(torsion computed using these two
Hilbert structures). This follows from the transformation
formula \ref{main properties of torsion for chain complexes}.
\end{remark}

The main result of this section is:
\begin{theorem} \label{torsion = filtered torsion}
We get with respect to the data \ref{data for chain complexes} and
the conventions above
$$\rho(C) ~= \rho_{\fil}(F_{\ast}C).$$
\end{theorem}

\noindent
It will follow from the next three lemmas.
\begin{lemma} \label{first reduction}
It suffices to treat the following special case:
$C$ itself
is a finite Hilbert cochain complex with the simple structure
represented by
\mbox{$\id : C \longrightarrow C$}
and the Hilbert structures on $F_p$ and $F_p/F_{p+1}$,
are obtained by the
given one on $C$ by taking Hilbert sub- and quotient structures.
\end{lemma}
\proof One easily constructs a finite Hilbert cochain complex
$D$ with
finite cofiltration $F_{\ast} D$ together with a chain homotopy
equivalence
\mbox{$f : D \longrightarrow C$} with the following property:
$f$ induces chain homotopy equivalences
\mbox{$F_pf: F_p D \longrightarrow F_p C$}
for all $p$ such that the given simple structure on
\mbox{$F_pC/F_{p+1}C$}
is represented by
\mbox{$F_pf/F_{p+1}f: F_p D/F_{p+1}D \longrightarrow F_p C/F_{p+1}C$}.
Equip $H(D)$ with the Hilbert structure for which
\mbox{$H(f) : H(D) \longrightarrow H(C)$} becomes an isometry.
We have by assumption and construction for all $p$ that
$$\rho(F_{p+1}D \longrightarrow F_pD \longrightarrow
F_pD/F_{p+1}D) ~ = ~ 0 ~ = ~
\rho(F_{p+1}C \longrightarrow F_pC \longrightarrow
F_pC/F_{p+1}C).$$
Hence we conclude from
\ref{main properties of torsion for chain complexes}
$$\rho(C) ~ = ~ \rho(D).$$
The maps $F_pf$ induce isometries between the
$E^{p + r\ast,q-(r-1)\ast}_r$ for $r \ge 1$ and
$H(E^{p + r\ast,q-(r-1)\ast}_r)$ for $r \ge 0$ associated
to $F_{\ast}D$ and $F_{\ast}C$. Now one easily verifies
$$\rho_{\fil}(F_{\ast}C) ~ = ~ \rho_{\fil}(F_{\ast}D).$$
This finishes the proof of Lemma \ref{first reduction}.
\qed

From now on we will only consider the special case described in
Lemma \ref{first reduction}. The proof of the next lemma is similiar to
the proof in \cite[Theorem 2.2.]{Munkholm(1981)}

\begin{lemma} \label{main step}
$$
\rho(C) ~ =~  \sum_{p} \rho(F_p/F_{p+1}) + \sum_{r \ge 1}
\sum_{p,q} (-1)^{p+q} \cdot [[\gamma^{p,q}_r]]
 - \sum_{p,q} (-1)^{p+q} \cdot [[\psi^{p,q}]]$$
\end{lemma}
\proof We do induction over $l$ with \mbox{$F_lC = \{0\}$}.
The begin of induction  $l=0$ is trivial, the induction step from $l-1$
to $l \ge 1$ done as follows.\par

Let $\overline{C}$ be $F_1 = F_1C$ with the cofiltration
$$\overline{F}_p ~ = ~ F_p\overline{C} ~ :=~  F_{p+1}C.$$
We will denote the various data coming from
the spectral sequence associated to $\overline{C}$ as the ones for
$C$ decorated with an additional bar.
Let $L\!H\!S$ be the acyclic finite Hilbert cochain complex given
by the long cohomology sequence associated to
\mbox{$0 \longrightarrow \overline{C} \longrightarrow C \longrightarrow
C/\overline{C} \longrightarrow 0$}.
It induces for $n \ge 0$
an acyclic  finite Hilbert cochain complex $D(n)$
concentrated in dimensions
$0$, $1$, $2$ and $3$ by
$$ H^n(C)/F^{1,n-1} \longrightarrow H^n(C/\overline{C}) \longrightarrow
H^{n+1}(\overline{C}) \longrightarrow F^{1,n} $$
where we equip $F^{1,n-1} = \im(H^n(\overline{C})
\longrightarrow H^n(C))$
respectively $H^{n+1}(C)/F^{1,n}$ with the Hilbert
substructure respectively quotient structure.
Define an acyclic
 Hilbert subcochain complex $D(n,r)$ for $r \ge 1$ of $D(n)$ by
$$H^n(C)/F^{1,n-1} \longrightarrow Z^{0,n}_r \longrightarrow
\overline{F}^{r-1,n+2-r} \longrightarrow F^{r,n+1-r} $$
Notice that $D(n,1) = D(n)$.
The following diagram  of acyclic finite Hilbert cochain complexes
concentrated in dimensions $1$, $2$ and $3$ commutes:
\begin{center}
$\begin{CD}
Z_r^{0,n}/Z_{r+1}^{0,n}  @>>>
\overline{F}^{r-1,n+2-r}/\overline{F}^{r,n+1-r}  @>>>
F^{r,n+1-r}/F^{r+1,n-r}
\\
@V\gamma_r^{0,n} VV  @V \overline{\psi}^{r-1,n+2-r} VV
@V \psi^{r,n+1-r} VV
\\
B_{r+1}^{r,n+1-r}/B_r^{r,n+1-r}  @>>>
\overline{Z}_{\infty}^{r-1,n+2-r}/\overline{B}_{\infty}^{r-1,n+2-r}
@>>>  Z_{\infty}^{r,n+1-r}/B_{\infty}^{r,n+1-r}
\end{CD}$
\end{center}
where the upper row is the quotient $D(n,r)/D(n,r+1)$ and
the lower row
is induced by the obvious inclusions and projections if one
takes the following identities into account
\begin{eqnarray*}
\overline{B}_{\infty}^{r-1,n+2-r} & = & B_r^{r,n+1-r};
\\
B_{r+1}^{r,n+1-r} & = & B_{\infty}^{r,n+1-r};
\\
\overline{Z}_{\infty}^{r-1,n+2-r} & = & Z_{\infty}^{r,n+1-r}.
\end{eqnarray*}
We conclude from the sum formula and transformation formula
\ref{main properties of torsion for chain complexes}
\begin{eqnarray*}
\rho(C) & = & \rho(\overline{C}) + \rho(C/\overline{C}) + \rho(L\!H\!S);
\\
\rho(L\!H\!S) & = & \sum_{n \ge -1} (-1)^{n+1} \rho(D(n));
\\
\rho(D(n)) & = &  \rho(D(n,l)) + \sum_{r = 1}^{l-1} \rho(D(n,r)/D(n,r+1));
\\
\rho(D(n,r)/D(n,r+1)) & = &
-[[\gamma_r^{0,n}]] + [[\overline{\psi}^{r-1,n+2-r}]] -
[[\psi^{r,n+1-r}]] \hspace{10mm} \mbox{ for } 1 \le r \le l.
\end{eqnarray*}
Since $D(n,l)$ is concentrated in dimensions $0$ and $1$
and its zero-th differential is $\psi^{0,n}$, we get
$$\rho(D(n,l)) ~ = ~ [[\psi^{0,n}]].$$
We compute using  the fact
$$\overline{\gamma}_r^{p,q} = \gamma_r^{p+1,q-1}\hspace{10mm}
\mbox{ for } p \ge 1$$
and the induction hypothesis applied to $\overline{C}$.
\begin{eqnarray*}
\rho(C) & = & \rho(\overline{C}) + \rho(F_0/F_1) +
\\
& & +\sum_{n \ge -1} (-1)^{n+1} \cdot  \left( [[\psi^{0,n}]] +
\sum_{r=1}^{l-1}
 - [[\gamma_r^{0,n}]] + [[\overline{\psi}^{r-1,n+2-r}]] -
[[\psi^{r,n+1-r}]]   \right)
\\[1mm]
& = & \sum_{p \ge 1} \rho(F_p/F_{p+1})  +
\sum_{r \ge 1} \sum_{p \ge 1,q} (-1)^{p+q-1} [[\gamma_r^{p,q-1}]]
- \sum_{p,q} (-1)^{p+q}\cdot
 [[\overline{\psi}^{p,q}]]
+ \rho(F_0/F_1)
\\
& & +\sum_{n \ge -1} (-1)^{n+1} \cdot  \left( [[\psi^{0,n}]] +
\sum_{r=1}^{l-1}
 - [[\gamma_r^{0,n}]] + [[\overline{\psi}^{r-1,n+2-r}]] -
[[\psi^{r,n+1-r}]]   \right)
\\[1mm]
& = & \sum_{p} \rho(F_p/F_{p+1}) + \sum_{r \ge 1}
\sum_{p,q} (-1)^{p+q} \cdot [[\gamma^{p,q}_r]]
 - \sum_{p,q} (-1)^{p+q} \cdot [[\psi^{p,q}]].
\end{eqnarray*}
This finishes the proof of Lemma \ref{main step}. \qed

\begin{lemma} \label{final step}
$$
\rho_{\fil}(F_{\ast}C) ~ = ~  \sum_{p} \rho(F_p/F_{p+1}) + \sum_{r \ge 1}
\sum_{p,q} (-1)^{p+q} \cdot [[\gamma^{p,q}_r]]
 - \sum_{p,q} (-1)^{p+q} \cdot [[\psi^{p,q}]]$$
\end{lemma}
\proof
We get from \ref{special Hilbert structure on cohomology}
$$\rho(E_r^{p+r\ast,q-(r-1)\ast}) ~ = ~
\sum_{k} (-1)^k \cdot [[\gamma_r^{p+rk,q-(r-1)k}]].$$
This implies
\begin{eqnarray*}
\sum_{r \ge 1}
\sum_{p,q} (-1)^{p+q} \cdot [[\gamma^{p,q}_r]] & = &
\sum_{r \ge 1}
\sum_{p = 0}^{r-1} \sum_q  (-1)^{p+q} \cdot
\sum_{k} (-1)^k \cdot [[\gamma_r^{p+rk,q-(r-1)k}]]
\\ & = & \sum_{r \ge 1}
\sum_{p = 0}^{r-1} \sum_q  (-1)^{p+q} \cdot
\rho(E_r^{p+r\ast,q-(r-1)\ast})
\end{eqnarray*}
This finishes the proof of Lemma \ref{final step} and hence of
Theorem \ref{torsion = filtered torsion}. \qed


\refstepcounter{section}
\typeout{----------- section 5 --------------------}

\tit{The fibration formula for Milnor torsion}
\label{The fibration formula for Milnor torsion}

For the remainder of this section we suppose that we have
a fibration \mbox{$f : E \longrightarrow B$}
 so that $B$ is a finite connected $CW$-complex,
the homotopy fiber has the homotopy type of a
finite $CW$-complex, the class
 \mbox{$\theta_f \in H^1(B;\Wh(E))$} defined in
\ref{cohomology class of fibration}
is trivial and $E$ is  locally path-connected.
\brn \label{fibration data}
Moreover, we assume that we are also given the following data:
\begin{enumerate}

\item A local coefficient system $V$ on $E$ with a unimodular Hilbert
structure;

\item a Hilbert structure on $H^{\ast}_{\sing}(E;V)$;

\item a coherent choice of simple structures on the fibers;

\item a unimodular Hilbert structure on
$\sum_{q} (-1)^q \cdot H^q_{\sing}(F;V)$
for the local coefficient systems
$H^q_{\sing}(F;V)$ over $B$;

\item a Hilbert structure on $H^{\ast}_{\cell}(B;H_{\sing}^q(F;V))$. \qed
\end{enumerate}
\ern

We mention that a choice of a unimodular Hilbert structure on
$\sum_{q} (-1)^q \cdot H^q_{\sing}(F;V)$  is possible because of

\begin{lemma} \label{unimodularity of sum H(F;V)}
If $f : E \longrightarrow B$ is a
fibration as desribed above, then
\mbox{$\sum_{q} (-1)^q \cdot H^q_{\sing}(F;V)$} is unimodular.
\end{lemma}
\proof We get from Lemma \ref{comparing Whitehead and Milnor torsion}
for a loop $w$ in $B$ with base point $b$ if
\mbox{$t_w : F_b \longrightarrow F_b$} is
given by the fiber transport along $w$:
$$\sum_{q} (-1)^q \cdot [[(t_w)^{\ast} : H^q_{\sing}(F_b;V)
\longrightarrow H^q_{\sing}(F_b;V)]]
~ = \rho(F_b;V) - \rho(F_b;V) + \Phi(V)(\theta_f(w)) ~ = ~ 0$$
where we think of $\theta_f$ as homomorphism
\mbox{$\pi_1(B,b) \longrightarrow \Wh(E)$}
and $\Phi(V)$ was introduced
in \ref{homomorphism from Wh to R}.
\qed\par

Notice that
\mbox{$\sum_{q} (-1)^q \cdot H^q_{\sing}(F;V)$}
is unimodular but not necessarily each \mbox{$H^q_{\sing}(F;V)$}.
\brn \label{explanation of torsion invariants}
We will consider the following real numbers:
\begin{enumerate}

\item $\rho(E;V)$\\
We get by construction \ref{simple structure on E}
a simple structure on $E$ from the data \ref{fibration data}.3
if we specify a cellular base point
system on $B$. However, by
\ref{change of simple structure under change of base point system}
the choice of cellular
base point system does not matter. Now we use
this simple structure and
data \ref{fibration data}.1 and \ref{fibration data}.2
to define the Milnor torsion
$\rho(E;V)$ according to
Definition \ref{Milnor torsion}. Notice that it depends only on
\ref{fibration data}.1,
\ref{fibration data}.2 and \ref{fibration data}.3.
This is the
invariant we want to compute;

\item $\rho(F;V)$\\
Choose $b \in B$. Then we get
from Definition \ref{Milnor torsion}
applied to data \ref{fibration data}.1,
\ref{fibration data}.3 and \ref{fibration data}.4 the Milnor torsion
$\rho(F_b;V)$. Here and elsewhere we supress in the notation
that we view $V$ as a local coefficient system over $F_b$ by the
inclusion of $F_b$ into $E$. We get from
Lemma \ref{comparing Whitehead and Milnor torsion} that
$\rho(F_b;V)$ is independent of $b$ since $B$ is path connected.
We abbreviate
\mbox{$\rho(F;V) = \rho(F_b;V)$}. This number depends only on
\ref{fibration data}.1, \ref{fibration data}.3 and
\ref{fibration data}.4;

\item $\rho(B;\sum_{q} (-1)^q \cdot H^q_{\sing}(F;V))$\\
This is the Milnor torsion and depends on the data
\ref{fibration data}.4 and \ref{fibration data}.5;

\item $\rho_{\sing}^{\Serre}(f)$\\
We have the Leray-Serre spectral sequence
for singular cohomology associated
to the fibration $f : E \longrightarrow B$. Namely,
the skeletal filtration of $B$ induces
a filtration $E_p = f^{-1}(B_p)$.  It yields a cofiltration
on $C^{\ast}_{\sing}(E;V)$ by putting
$$F_pC^{\ast}_{\sing}(E;V) ~ = ~ C^{\ast}_{\sing}(E,E_{p-1};V).$$
Recall that $C^{\ast}_{\sing}(E,E_{p-1};V)$ is
the kernel of the canonical map
given by the inclusion
\mbox{$C^{\ast}_{\sing}(E;V) \longrightarrow
C^{\ast}_{\sing}(E_{p-1};V)$}.
We will later recall in \ref{identification of E_2-term}
the isomorphism which identifies the $E_2$-term
of the associated spectral sequence
$$U_2^{p,q} : H^p_{\cell}(B;H^q_{\sing}(F;V))
\longrightarrow E^{p,q}_2.$$
Equip $E_2^{p,q}$ with the Hilbert structure so that $U_2^{p,q}$ becomes an
 isometry.
Using definition 
\ref{torsion of a filtered chain complex} put
$$\rho_{\sing}^{\Serre}(f) := 
\rho_{\fil}^{\ge 2}(F_{\ast}C^{\ast}_{\sing}(X;V)) + \sum_{p,q}(-1)^{p+q}[[U_2^{p,q}]].$$
Using remark \ref{def of rho ge 2}, $\rho_{\sing}^{\Serre}(f)$ depends only on the data
\ref{fibration data}.2 and \ref{fibration data}.5.
\end{enumerate}\ern

The main result of this section is

\begin{theorem} \label{fibration formula}
If $F$ is a fibration as described above with data
\ref{fibration data} then
$$\rho(E;V) ~ = ~ \chi(B) \cdot \rho(F;V) +
\rho(B;\sum_{q} (-1)^q \cdot H^q_{\sing}(F;V)) +
\rho_{\sing}^{\Serre}(f).$$
\end{theorem}
\proof

Fix a cellular base point system
\mbox{$\{b_c \mid c \in I\}$} for $B$ and an
orientation for each cell $c \in I$. The homotopy equivalence of
Lemma \ref{simple structure on (E_n,E_(n-1)} together with the
suspension cochain homotopy equivalence and the obvious
identification of $F_p/F_{p+1}$ with $C^{\ast}_{\sing}(E_p,E_{p-1};V)$
yields a cochain homotopy equivalence unique up to homotopy
$$U_0^{p,\ast} : \oplus_{c \in I_p} \Sigma^p C^{\ast}_{\sing}(F_{b_c};V)
\longrightarrow F_p/F_{p+1}.$$
We have defined an isomorphism in \ref{identification of modules}
$$\mu^{p,q} : C^p_{\cell}(B;H^q_{\sing}(F;V)) \longrightarrow
\oplus_{c \in I_p} H^q_{\sing}(F_{b_c};V).$$
Recall that we equip
\mbox{$C^p_{\cell}(B;H^q_{\sing}(F;V))$}
with the Hilbert structure for which $\mu^{p,q}$ becomes an isometry.
Define an isomorphism by the composition
$$U^{p,q}_1 : C^p_{\cell}(B;H^q_{\sing}(F;V)) \xrightarrow{\mu^{p,q}}
\oplus_{c \in I_p} H^q_{\sing}(F_{b_c};V) =
\oplus_{c \in I_p} H^{p+q}(\Sigma^p C^{\ast}_{\sing}(F_{b_c};V))$$
$$\xrightarrow{H^{p+q}(U_0^{p,\ast})} H^{p+q}(F_p/F_{p+1}) =
H^0(E_0^{p,q+\ast}) \xrightarrow{\phi^{p,q}_0} E_1^{p,q}$$
This isomorphism is independent of the choice of base point system and
orientation of the cells and is compatible with the differentials on
\mbox{$C^{\ast}_{\cell}(B;H^q_{\sing}(F;V))$} and $E_1^{\ast,q}$.
Hence we get a canonical isomorphism
\begin{eqnarray}
U^{p,q}_2 : H^p_{\cell}(B;H^q_{\sing}(F;V))
\xrightarrow{H^p(U_1^{\ast,q})}
& H^p(E_1^{\ast,q}) = H^0(E_1^{p+\ast,q}) &
\xrightarrow{\phi_1^{p,q}} E^{p,q}_2.
\label{identification of E_2-term}
\end{eqnarray}
The simple structure on $F_{b_c}$ induces one on
\mbox{$C^{\ast}_{\sing}(F_{b_c};V)$}. Put on $F_p/F_{p+1}$
the simple structure for which the torsion
\mbox{$t(U_0^{p,\ast} :
\oplus_{c \in I_p} \Sigma^p C^{\ast}_{\sing}(F_{b_c};V)
\longrightarrow F_p/F_{p+1})$}  defined in
\ref{torsion for simple chain structures} vanishes. Equip $F_p$
inductively with the simple structure for which
\mbox{$\rho(F_{p+1} \longrightarrow F_p \longrightarrow F_p/F_{p+1})$}
defined in \ref{torsion of exact sequence} vanishes.
Then this simple structure on \mbox{$F_0 = C^{\ast}_{\sing}(E;V)$}
agrees with the one we get from the simple structure on $E$
which we have defined
in \ref{simple structure on E} by an inductive procedure
over the $E_n$-s  and Lemma \ref{simple structure on (E_n,E_(n-1)}.
We conclude from
Theorem \ref{torsion = filtered torsion}:
$$\rho(E;V) ~ = ~ \rho_{\fil}(F_{\ast}C^{\ast}_{\sing}(E;V)).$$
Because $\phi_r^{p,q}$ are isometries, 
we get from the transformation formula
\ref{main properties of torsion for chain complexes}:
\begin{eqnarray*}
\sum_{p}\rho(F_p/F_{p+1}) & = &
\sum_{p} \sum_{c \in I_p}
\rho(\Sigma^p C^{\ast}_{\sing}(F_{b_c};V)) +
\sum_p \sum_q (-1)^q \cdot [[H^q(U_0^{p,\ast})]]
\\
&  = & \sum_p (-1)^p \cdot \sum_{c \in I_p} \rho(F;V) +
\sum_{p,q} (-1)^{p+q} \cdot [[H^{p+q}(U_0^{p,\ast})]]
\\
& = & \chi(B) \cdot \rho(F;V) + \sum_{p,q} (-1)^{p+q} \cdot [[U_1^{p,q}]]\\
\sum_{q} (-1)^q \cdot \rho(E_1^{\ast,q}) & = &
\sum_q(-1)^q\cdot \rho(C^{\ast}_{\cell}(B;H^q_{\sing}(F;V)) -
\sum_{q} (-1)^q \cdot \sum_p (-1)^p \cdot [[U_1^{p,q}]] \\
& & + \sum_q (-1)^q \cdot \sum_p (-1)^p \cdot [[H^p(U_1^{\ast,q})]]
\\
 & = & \rho(B;\sum_q (-1)^q \cdot H^q(F;V)) -
\sum_{p,q} (-1)^{p+q} \cdot [[U_1^{p,q}]] + \sum_{p,q} (-1)^{p+q} \cdot [[U_2^{p,q}]]
\end{eqnarray*}
We get from
\ref{explanation of torsion invariants}.4 and
Definition \ref{torsion of a filtered chain complex}:
\begin{eqnarray*}
\rho_{\fil}(F_{\ast}C^{\ast}_{\sing}(E;V)) & = &
\sum_{p} \rho(F_p/F_{p+1}) + \sum_{q} (-1)^q \cdot
\rho(E_1^{\ast,q})
\\
& &  +  \rho_{\sing}^{\Serre}(f) - \sum_{p,q} (-1)^{p+q} \cdot [[U_2^{p,q}]].
\end{eqnarray*}
Now Theorem \ref{fibration formula} follows. \qed

\begin{remark} \label{simple structure on B suffices again} \em
 One can extend Theorem \ref{fibration formula} to the case where
$B$ is not necessarily a finite $CW$-complex
but carries a simple structure
and one requires in data \ref{fibration data}.5 a Hilbert structure on
$H^{\ast}_{\sing}(B;H^q(F;V))$ instead of
$H^{\ast}_{\cell}(B;H^q(F;V))$. Then we still get a simple
structure on $E$ by
Remark \ref{simple structure on B suffices}. It remains to define
$\rho^{\Serre}_{\sing}(f)$. Choose an arbitrary homotopy equivalence
\mbox{$h : Y \longrightarrow B$} for a finite $CW$-complex $Y$. Let
\mbox{$\overline{f} : h^{\ast}E \longrightarrow Y$}
be the pull back fibration
and \mbox{$\overline{h} : h^{\ast}E \longrightarrow E$}
 be the canonical fiber
homotopy equivalence. Equip
\mbox{$H^{\ast}_{\sing}(h^{\ast}E;\overline{h}^{\ast}V)$}
and
\mbox{$H^{\ast}_{\cell}(Y;h^{\ast}H^q_{\sing}(F;V))$}
with the Hilbert structures for which the following
isomorphisms become isometries (see \ref{comparing singular and cellular})
$$H^{\ast}_{\sing}(h^{\ast}E;\overline{h}^{\ast}V)
\xrightarrow{(\overline{h}^{\ast})^{-1}}
H^{\ast}_{\sing}(E;V)$$
and
$$H^{\ast}_{\cell}(Y;h^{\ast}H^q_{\sing}(F;V)) \longrightarrow
H^{\ast}_{\sing}(Y;h^{\ast}H^q_{\sing}(F;V)) \xrightarrow{(h^{\ast})^{-1}}
H^{\ast}_{\sing}(B;H^q_{\sing}(F;V)).$$
Then $\rho^{\Serre}_{\sing}(\overline{f})$ is defined and we put
$$\rho^{\Serre}_{\sing}(f) ~ = ~ \rho^{\Serre}_{\sing}(\overline{f}).$$
This is independent of the choice of $h$ and Theorem
\ref{fibration formula} remains true.

\end{remark}



\refstepcounter{section}
\typeout{----------- section 6 --------------------}

\tit{The fiber bundle formula for analytic torsion}
\label{The fiber bundle formula for analytic torsion}

In this section we give the proof of
Theorem \ref{fiber bundle formula for analytic torsion}
and  Corollaries \ref{acyclic case},
\ref{spheres as fibers} and
\ref{spheres as base}.

We begin with the proof of Theorem
\ref{fiber bundle formula for analytic torsion}.
In the first step we
reduce the claim to the case where the Riemannian  metric on $V$
is unimodular. By assumption we only have to treat the case where
$\dim(E)$ is odd. We will vary the Riemannian metric on $V$,
but fix the Riemannian metrics on $E$ and $B$.
Denote by $V$ the flat bundle $V$ with
the given Riemannian metric and by $\widehat{V}$ the flat
vector bundle $V$ with some unimodular Riemannian metric.
Let $\rho^{\Serre}_{\dR}(f)$ and $\widehat{\rho^{\Serre}_{\dR}(f)}$
be the two values of the correction term for the choice of these two
Riemannian metrics on $V$. Notice that either $\Pf_B$ is zero or
$\dim(F_b)$ is odd.
We get from theorem \ref{variation of Riemanninan metrics}
and definition \ref{definition of rho^Serre as real number}
\begin{eqnarray*}
\rho_{\an}(E;V) - \rho_{\an}(E;\widehat{V})
& = &
- \sum_p (-1)^p \cdot [[H^p_{\dR}(E;V)
\xrightarrow{\id} H^p_{\dR}(E;\widehat{V})]];
\\
\makebox[20mm]{$\int_B \rho_{\an}(F_b;V) \cdot \Pf_B
- \int_B \rho_{\an}(F_b;\widehat{V}) \cdot \Pf_B$} & &
\\
 & = &
- \int_B \left(\sum_q (-1)^q \cdot [[H^q_{\dR}(F_b;V) \xrightarrow{\id}
H^q_{\dR}(F_b;\widehat{V})]]\right) \cdot \Pf_B;
\\
\makebox[20mm]{$\rho_{\an}(B;H^q_{\dR}(F;V))-
\rho_{\an}(B;H^q_{\dR}(F;\widehat{V}))$} & &
\\
& = & - \sum_p (-1)^p \cdot [[H^p_{\dR}(B;H^q_{\dR}(F;V)) \xrightarrow{\id}
H^p_{\dR}(B;H^q_{\dR}(F;\widehat{V}))]]
\\
& &
+ \int_B [[H^q_{\dR}(F_b;V) \xrightarrow{\id} H^q_{\dR}(F_b;\widehat{V})]]
\cdot \Pf_B;
\\
\rho^{\Serre}_{\dR}(f) - \widehat{\rho^{\Serre}_{\dR}(f)}
& = &
\sum_{p,q} (-1)^{p+q} \cdot [[H^p_{\dR}(B;H^q_{\dR}(F;V))
\xrightarrow{\id}
H^p_{\dR}(B;H^q_{\dR}(F;\widehat{V}))]]
\\
& & -
\sum_p (-1)^p \cdot [[H^p_{\dR}(E;V) \xrightarrow{\id}
H^p_{\dR}(E;\widehat{V})]]
\end{eqnarray*}

We conclude from the equations above:
\begin{eqnarray*}
 & & \rho_{\an}(E;V) - \left(\int_B \rho_{\an}(F_b;V) \cdot \Pf_B
 ~ + ~ \sum_{q} (-1)^q \cdot \rho_{\an}(B;H^q_{\dR}(F;V))
 ~ + \rho^{\Serre}_{\dR}(f)\right)
\\
 & = & \rho_{\an}(E;\widehat{V}) -
\left(\int_B \rho_{\an}(F_b;\widehat{V}) \cdot \Pf_B
 ~ + ~ \sum_{q} (-1)^q \cdot \rho_{\an}(B;H^q_{\dR}(F;\widehat{V}))
 ~ + \widehat{\rho^{\Serre}_{\dR}(f)}\right).
\end{eqnarray*}
Hence we can assume in the sequel without loss of generality
that the Riemannian metric on $V$ is unimodular.\par

Given a closed smooth manifold $M$ and a flat vector bundle
$W$ over $M$ with Riemannian metrics on $M$ and $W$,
the {\em harmonic Hilbert structure } on the singular cohomology
\mbox{$H^{\ast}_{\sing}(M;W)$} is given by the deRham isomorphism
$$H^{\ast}_{\dR}(M;W) \longrightarrow H^{\ast}_{\sing}(M;W)$$
and the harmonic Hilbert structure on $H^{\ast}_{\dR}(M;W)$
coming from the Hodge decomposition and the Hilbert space structure
on the space of harmonic forms as explained
in the introduction. If we have choosen a
smooth triangulation on $M$, the harmonic Hilbert structure on
\mbox{$H^{\ast}_{\cell}(M;W)$} is induced by the canonical isomorphism
$$H^{\ast}_{\cell}(M;W) \longrightarrow H^{\ast}_{\sing}(M;W)$$
induced from \ref{comparing singular and cellular}.
Denote by $\rho_{\topological}(M;W)$ the Milnor torsion of
$M$ defined in \ref{Milnor torsion}
with respect to the  simple structure represented by a
smooth triangulation and the harmonic Hilbert structure on
\mbox{$H^{\ast}_{\sing}(M;V)$}.\par

Denote by $\overline{H^q_{\sing}(F;V)}$ the flat bundle over $B$
with fibers ${H^q_{\sing}(F_b;V)}$ equipped with
some Riemannian metric such that the induced Hilbert structure
on \mbox{$\sum_q (-1)^q \cdot \overline{H^q_{\sing}(F;V)}$} is
unimodular. Recall that this can be done by Lemma
\ref{unimodularity of sum H(F;V)}. Equip $B$ and $F_b$ with the
$CW$-structure given by some smooth triangulation. Let
\mbox{$\overline{\rho(F_b;V)}$} be the Milnor torsion
of $F_b$ with respect to  $\overline{H^q_{\sing}(F;V)}$ (see Definition
\ref{Milnor torsion}). We have seen in
\ref{explanation of torsion invariants}.2 that it is
independent of $b \in B$ and abbreviate it by
\mbox{$\overline{\rho(F;V)}$}.
Let \mbox{$\overline{\rho^{\Serre}_{\sing}(f)}$} be the correction
term of \ref{explanation of torsion invariants} with respect to
the harmonic Hilbert structures on $H^{\ast}_{\sing}(E;V)$ and
\mbox{$H^p_{\sing}(B;\overline{H^q_{\sing}(F;V)})$}.
We conclude from
Lemma \ref{smooth bundles and simple structures} and
Theorem \ref{fibration formula}:
\begin{eqnarray}
\rho_{\topological} (E;V) & = &
\chi(B) \cdot \overline{\rho(F;V)} +
\rho(B;\sum_{q} (-1)^q \cdot \overline{H^q_{\sing}(F;V)}) +
\overline{\rho_{\sing}^{\Serre}(f)}.
\label{topological formula}
\end{eqnarray}
Next we want to show:
\begin{eqnarray}
\rho_{\an}(E;V) & = & \rho_{\topological}(E;V);
\label{top = an formula 1}
\\
\rho_{\an}(F_b;V) & = & \rho_{\topological}(F_b;V);
\label{top = an formula 2}
\\
\sum_{q} (-1)^q \cdot \rho_{\an}(B;\overline{H^q_{\sing}(F;V)})
& = &
\rho(B;\sum_{q} (-1)^q \cdot \overline{H^q_{\sing}(F;V)}).
\label{top = an formula 3}
\end{eqnarray}

The first two equations follow directly from Theorem
\ref{equality of analytic and topological torsion}. For the
proof of the last equation we must take a closer look at the
result of Bismut and Zhang
\cite[Theorem 0.2]{Bismut-Zhang (1992)}. The problem is
that the Riemannian metric on
\mbox{$\sum_q (-1)^q \cdot
 \overline{H^q_{\sing}(F;V)}$} is unimodular
but it is not true in
general that each of the flat bundles \mbox{$H^q_{\sing}(F;V)$} is
unimodular or, equivalently, admits some unimodular Riemannian
metric.\par

Choose a Morse function $f$ on $B$. Let $X$ be its gradient vector
field with respect to some Riemannian metric on $B$
which can be different from the given Riemannian metric such that $X$
satisfies the Smale transversality conditions. We obtain a
$CW$-structure on $B$, denoted by $B^{\prime}$, and a cellular
base point system \mbox{$\{ b_c \mid c \in I\}$} given by
the critical points of $f$. The simple structure represented by
\mbox{$\id : B^{\prime} \longrightarrow B$} is the same as the one
represented by any smooth triangulation. Let
\mbox{$C^{\ast}_{\cell}(B^{\prime},
\{b_c\};\overline{H^q_{\sing}(F;V)})$}
be the cellular cochain complex with the Hilbert structure
defined in \ref{identification of modules}. We equip \mbox{$H^p(C^{\ast}_{\cell}(B^{\prime},\{b_c\};
\overline{H^q_{\sing}(F;V)}))$} with the harmonic Hilbert structure.
We get by  inspecting Definition \ref{Milnor torsion}
$$\rho(B,\sum_q (-1)^q \cdot \overline{H^q_{\sing}(F;V)}) ~ = ~
\sum_{q} (-1)^q \cdot \rho(C^{\ast}_{\cell}(B^{\prime},
\{b_c\};\overline{H^q_{\sing}(F;V)})).$$
But \mbox{$C^{\ast}_{\cell}(B^{\prime},
\{b_c\};\overline{H^q_{\sing}(F;V)})$} is isometrically isomorphic
to the Thom-Smale complex. Let  $\vol_{\harm}$
be the volume form associated to the harmonic Hilbert structure.
Hence we get from the definitions (cf.
\cite[Remark 2.3 and Remark 1.10]{Bismut-Zhang (1992)})
\begin{eqnarray*}
\rho_{\an}(B;\overline{H^q_{\sing}(F;V)}) & = &
- \ln\left(
||\vol_{\harm}||^{RS}_{\det H^{\bullet}_{\dR}
(B;\overline{H^q_{\sing}(F;V)})}\right);
\\
\rho(C^{\ast}_{\cell}(B^{\prime},
\{b_c\};\overline{H^q_{\sing}(F;V)}) & = &
 - \ln\left(||\vol_{\harm}||^{{\cal M},X}_{\det H^{\bullet}_{\cell}
 (B;\overline{H^q_{\sing}(F;V)})}\right);
\end{eqnarray*}
where the terms on the right side are the invariants appearing in
\cite{Bismut-Zhang (1992)}. We get from
\cite[Theorem 0.2]{Bismut-Zhang (1992)})
$$\rho_{\an}(B;\overline{H^q_{\sing}(F;V)}) -
\rho(C^{\ast}_{\cell}(B^{\prime},
\{b_c\};\overline{H^q_{\sing}(F;V)})) ~ = ~
\int_B \theta(\overline{H^q_{\sing}(F;V)}) \cdot X^{\ast} \psi(B)$$
for a certain  $(\dim(B) -1)$-form $\psi(B)$. We conclude
\begin{eqnarray*}
 & & \sum_q (-1)^q \cdot \rho_{\an}(B;\overline{H^q_{\sing}(F;V)}) -
 \rho(B,\sum_q (-1)^q \cdot \overline{H^q_{\sing}(F;V)}) \\
& = & \sum_q (-1)^q\cdot \rho_{\an}(B;\overline{H^q_{\sing}(F;V)}) -
\sum_{q} (-1)^q \cdot \rho(C^{\ast}_{\cell}(B^{\prime},
\{b_c\};\overline{H^q_{\sing}(F;V)}))
\\
& = & \int_B \sum_q (-1)^q \cdot \theta(\overline{H^q_{\sing}(F;V)})
\cdot X^{\ast} \psi(B).
\end{eqnarray*}
Now one easily checks that
\mbox{$\sum_q (-1)^q \cdot \theta(\overline{H^q_{\sing}(F;V)})$}
vanishes as the Riemannian metric on
\mbox{$\sum_q (-1)^q \cdot \overline{H^q_{\sing}(F;V)})$} is unimodular.
This finishes the proof of \ref{top = an formula 3}.\par

Denote by $H^q_{\dR}(F;V)$ and $H^q_{\sing}(F;V)$ the flat
bundles with the harmonic Riemannian metrics (in contrast to
$\overline{H^q_{\sing}(F;V)}$). We get from the
transformation formula
\ref{main properties of torsion for chain complexes}
and \ref{top = an formula 2}
\begin{eqnarray} \label{Formel eins}
& & \chi(B) \cdot \overline{\rho(F;V)} - \int_B \rho_{\an}(F_b;V)
\cdot \Pf_B
\\
& = &
\int_B \left( \overline{\rho(F_b;V)} - \rho_{\topological}(F_b;V)\right)
\cdot \Pf_B \nonumber
\\
 & = &   - \int_B \sum_q (-1)^q \cdot
[[\overline{H^q_{\sing}(F_b;V)} \xrightarrow{\id}
H^q_{\sing}(F_b;V)]] \cdot \Pf_B. \nonumber
\end{eqnarray}
We conclude from \ref{top = an formula 3} and
Theorem \ref{variation of Riemanninan metrics}:
\begin{eqnarray}
\label{Formel Zwei}
 & & \rho(B;\sum_{q} (-1)^q \cdot \overline{H^q_{\sing}(F;V)})
- \sum_{q} (-1)^q \cdot \rho_{\an}(B;H^q_{\dR}(F;V))\\
& = &
\sum_{q} (-1)^q \cdot
\left(\rho_{\an}(B;\overline{H^q_{\sing}(F;V)}) -
\rho_{\an}(B;H^q_{\sing}(F;V))\right) \nonumber
\\
 & = & \sum_q (-1)^q \cdot \left( - \sum_p (-1)^p \cdot
[[H^p_{\dR}(B,\overline{H^q_{\sing}(F;V)})
\xrightarrow{\id} H^p_{\dR}(B,H^q_{\sing}(F;V))]] \right. \nonumber
\\
& & \left.+
\int_B [[\overline{H^q_{\sing}(F;V)} \xrightarrow{\id} H^q_{\sing}(F;V)]]
\cdot \Pf_B\right)\nonumber
\\
& = &
- \sum_{p,q} (-1)^{p+q} \cdot
[[H^p_{\dR}(B,\overline{H^q_{\sing}(F;V)}) \xrightarrow{\id}
H^p_{\dR}(B,H^q_{\sing}(F;V))]] \nonumber
\\
 & & +
\sum_q (-1)^q \cdot \int_B [[\overline{H^q_{\sing}(F_b;V)}
 \xrightarrow{\id} H^q_{\sing}(F_b;V)]] \cdot \Pf_B. \nonumber
\end{eqnarray}

The deRham isomorphism
\mbox{$H^n_{\dR}(E;V) \longrightarrow H^n_{\sing}(E;V)$}
is compatible with the two filtrations and the identifications
of the $E_2$-term in the Leray-Serre spectral sequence for
deRham and singular cohomology
\ref{identification of E_2term for deRham} and
\ref{identification of E_2-term}. It yields isomorphisms on $E_r^{p,q}$ for 
$r\ge 2$. This implies using remark \ref{def of rho ge 2}
\begin{equation}
\overline{\rho^{\Serre}_{\sing}(f)} -
\rho^{\Serre}_{\dR}(f) =
\sum_{p,q} (-1)^{p+q}
[[H^p_{\dR}(B,\overline{H^q_{\sing}(F;V)}) \xrightarrow{\id}
H^p_{\dR}(B,H^q_{\sing}(F;V))]].
\label{Formel Drei}
\end{equation}

Now Theorem \ref{fiber bundle formula for analytic torsion}
follows from \ref{topological formula}, \ref{top = an formula 1},
\ref{Formel eins}, \ref{Formel Zwei} and \ref{Formel Drei}. \qed

Next we prove Corollary \ref{acyclic case}. Recall that
\mbox{$H^q_{\dR}(F;V)$} vanishes by assumption. In particular the
$E_2$-term of the Leray-Serre spectral sequences for cohomology vanishes.
Hence Theorem \ref{fiber bundle formula for analytic torsion} implies
$$\rho_{\an}(E;V) ~ = ~ \int_B \rho_{\an}(F_b;V) \cdot \Pf_B.$$
If $\dim(B)$ is odd, the right hand side of the equation above and
$\chi(B) \cdot \rho_{\an}(F_b;V)$ vanish and the claim follows.
Suppose that $\dim(B)$ is even. Then $\dim(F_b)$ is even and the
Riemannian metric on $V$ is unimodular  by assumption or
$\dim(F_b)$ is odd. In both cases $\rho_{\an}(F_b;V)$ is independent
of $b \in B$ by Theorem \ref{variation of Riemanninan metrics} and the claim follows.
This finishes the proof of Corollary \ref{acyclic case}. \qed

Next we outline the proof of Corollary
\ref{spheres as fibers}. Analogously to
the first step in the proof of Theorem
\ref{fiber bundle formula for analytic torsion} we can show that we
can assume without loss of generality that the Riemannian metric
on $W$ is unimodular. Put $V = f^{\ast}W$.
There is a canonical isomorphism of local coefficient
systems over $B$
$$H^q_{\sing}(F;\zz) \otimes_{\zz} W \longrightarrow
H^q_{\sing}(F;V).$$
Recall that we assume that $E$ and $B$ are oriented.
(One can drop this assumption if one allows an additional
twist for $W$.) This implies that $\pi_1(B)$ acts orientation
preserving on the homology of the fiber. Hence we can
choose an isomorphism of local coefficient systems of
$\zz$-module
from $H^q_{\sing}(F;V)$ to the trivial coefficient system
with value $\zz$. Hence we obtain identifications
$$ H^q_{\sing}(F;V) ~ = ~ \left\{
\begin{array}{lr}
W  &  q =0,n \\ 0 & q \not= 0,n
\end{array} \right..$$
We conclude from Theorem \ref{fibration formula}
$$\rho_{\topological}(E;V) ~ = ~
\chi(B) \cdot \rho(F_b;V) +
\rho(B;\sum_q (-1)^q \cdot H^q_{\sing}(F;V)) +
\rho^{\Serre}_{\sing}(f)$$
where we use on $H^q_{\sing}(F_b;V)$  the unimodular Hilbert structure
given by the one on $W$ and the identification above and on
\mbox{$H^q_{\sing}(B;\sum_q (-1)^q \cdot H^q_{\sing}(F;V))$} the
harmonic one. We get
\begin{eqnarray*}
\rho(F;V) & = & 0;\\
\rho(B;\sum_q (-1)^q \cdot H^q_{\sing}(F;V))
& = & \chi(S^n) \cdot \rho_{\topological}(B;W);
\\
\rho^{\Serre}_{\sing}(f) & = & \rho_{\dR}^{\Serre}(f);
\\
\rho_{\an}(E;V) & = & \rho_{\topological}(E;V);
\\
\rho_{\an}(B;W) & = & \rho_{\topological}(B;W);
\end{eqnarray*}
since there is a  canonical isomorphism
$$C^{\ast}_{\sing}(F_b;\zz) \otimes_{\zz} W
\longrightarrow C^{\ast}_{\sing}(F_b;V)$$
and  $H^{\ast}_{\sing}(F_b;\zz)$ is
free as $\zz$-module, the deRham isomorphism is
compatible with the cofiltrations on the de Rham complex
and the singular cochain complex and with the identifications of
the $E_2$-terms of the associated Leray-Serre spectral sequences
\ref{identification of E_2term for deRham} and
\ref{identification of E_2-term}
and we have Theorem
\ref{equality of analytic and topological torsion}.
Hence it remains to show
$$\rho^{\Serre}_{\dR}(f) ~ = ~ \rho(G^{\ast}).$$
Notice that the $E^2$-term of the Leray-Serre spectral sequence is
trivial except for the $0$-th and $n$-th row.
Hence we can splice the
spectral sequence together to one exact sequence
$$\ldots \longrightarrow E^{p,n}_{n} \xrightarrow{d_{n+1}^{p,n}}
 E_{n}^{p+n+1,0}
\longrightarrow H^{p+n+1}(E;V) \longrightarrow E^{p+1,n}_{n}
\xrightarrow{d_{n+1}^{p+1,n}} \ldots $$
If one takes the identification of the $E_2$-term
\ref{identification of E_2term for deRham} and the obvious
identification
\mbox{$E^{p,q}_2 = E^{p,q}_n$}
into account,
one gets an exact sequence whose torsion is precisely
$\rho^{\Serre}_{\dR}(f)$. Moreover, it can be
identified with the Gysin sequence up to sign.
This finishes the proof of Corollary
\ref{spheres as fibers}.
The proof of Corollary \ref{spheres as base}
is similar.  \qed


\refstepcounter{section}
\typeout{----------- section: torsion + determinants --------------------}
\begin{appendix}
\section{Torsion and determinants}
\label{Torsion and determinants}

We start with a description of the determinant of a finite dimensioinal vector space.
Let $V$ be a $n$-dimension (real) vector space. Define its
{\em determinant} $\det(V)$ to be the $1$-dimensional vector space
$\wedge^nV$ given by the $n$-th exterior power. There are canonical
isomorphisms
\begin{eqnarray*}
\det(V \oplus W) & \xrightarrow{\cong} & \det(V) \otimes \det(W);
\\
\det(V^{\ast}) & \xrightarrow{\cong} & \det(V)^{\ast};
\\
\det(V)^{\ast} \otimes \det(W) & \xrightarrow{\cong} &
\hom_{\rr}(\det(V),\det(W));
\\
\det(U^{\ast}) \otimes \det(V) \otimes \det(V^{\ast}) \otimes \det(W)
& \xrightarrow{\cong} &
\det(U^{\ast}) \otimes \det(W);
\\
\det(V) & \longrightarrow (\det(V)^{\ast})^{\ast}.
\end{eqnarray*}
Given a homomorphisms \mbox{$f: V \longrightarrow W$}, we obtain a
well-defined element
$$\det(f) \in \det(V)^{\ast} \otimes \det(W).$$
Under the identifications above
\begin{eqnarray*}
\det(f \circ g) & = & \det(f) \otimes \det(g);
\\
\det\squarematrix{f}{h}{0}{g} & = & \det(f) \otimes \det(g).
\end{eqnarray*}
If $(V_p)_p$ is a graded vector space graded by nonnegative integers
\mbox{$p \ge 0 $} and almost all $V_p$ are zero, then we define
$$\det((V_p)_p) ~ := ~ \otimes_{p \ge 0} \det(V_p)^{(-1)^p}$$
where $\det(V)^1$ is $\det(V)$ and  $\det(V)^{-1}$ is defined by
$\det(V)^{\ast}$. If $(V_{p,q})_{p,q}$ is bigraded by non-negative integers
$p,q \ge 0$ and almost all $V_{p,q}$ are zero, then
we define
$$\det((V_{p,q})_{p,q}) ~ = ~ \det((W_n)_n)$$
where $(W_n)_n$ is the associated graded vector space with
\mbox{$W_n = \oplus_{p=0}^n V_{p,n-p}$}.

Now we will introduce algebraic torsion invariants in the language of determinants..
Given an acyclic finite cochain complex, we define
\begin{eqnarray}
\overline{\rho}(C) & \in & \det((C^p)_p)^{-1}
\label{torsion of an acyclic chain complex by determinants}
\end{eqnarray}
by the determinant of the isomorphism
\mbox{$c^{\ast} + \gamma^{\ast}: C^{\ev} \longrightarrow C^{\odd}$}.
For a homotopy equivalence
\mbox{$f: C \longrightarrow D$} of finite cochain complexes
we define
\begin{equation}
\overline{\rho}(f)  :=  \overline{\rho}(\cone(f))
 \in\det((C^p)_p)^{-1} \otimes \det((D^p)_p) \cong
\hom_{\rr}(\det(C^p)_p,\det(D^p)_p)
.
\label{definition of torsion by determinants}
\end{equation}
Given a finite cochain complex $C$, we use the cochain map
\mbox{$i: H(C) \longrightarrow C$} to define the isomorphism
\begin{eqnarray}
\overline{\rho}(C) & := & (\overline{\rho}(i))^{-1}
\hspace{10mm} : \det((C^p)_p) \longrightarrow \det((H^p(C))_p).
\label{definition of Milnor torsion by determinants}
\end{eqnarray}
The topological torsion $\overline{\rho}(X)$ of a finite CW-complex  $X$ 
is defined as the torsion of the
cellular cochain complex.

The definition of the torsion of a filtration/spectral sequence is a little bit
more elaborate.

Let  $F_{\ast}W$ be a fibration
of a finite-dimensional vector space $W$
$$\{0\} = F_{-1} \subset \ldots F_{p} \subset F_{p+1}
\subset \ldots \subset F_{n} = W. $$
This determines  an isomorphism
\begin{eqnarray}
\overline{\rho}^{\fil}(F_{\ast}W):
\otimes_{p \ge 0} \det(F_p/F_{p-1})) & \xrightarrow{\cong} & \det(W).
\label{torsion of a filtration}
\end{eqnarray}
which is defined inductively over $n$. The begin of induction
$n=1$ is given by the Milnor torsion of the acyclic
$2$-dimensional cochain complex
\mbox{$F_0 \longrightarrow W \longrightarrow W/F_0$}.
For the step $n-1$ to $n$ observe that 
\[ \overline{\rho}^{\fil}(F_*F_{n-1}: 
\bigoplus_{p=0}^{n-1} \det(F_p/F_{p-1})\xrightarrow{\cong}\det(F_{n-1} \]
is defined by induction, and
 we simply compose (after tensoring with $\det(W/F_{n-1})$) 
with the Milnor torsion of $F_{n-1}\to W\to W/F_{n-1}$.

In particular, suppose $f:E\to B$ is a fibration of closed manifolds as in the introduction.
This yields a filtration of $H^n_{\dR}(E;V)$ for $n \ge 0$.
Take the tensor product of the corresponding torsion isomorphisms
resp.~their inverses depending on the parity of $n$ to obtain
an isomorphism
$$\overline{\rho}^{\fil}: \det((F^{p,q}/F^{p+1,q-1})_{p,q})
\longrightarrow \det((H^n(E;V))_n).$$
In the de Rham spectral sequence, we can compute the torsions of
the cochain complex
\mbox{$(E_r^{p+r\ast,q-(r-1)\ast},d_r^{p+r\ast,q-(r-1)\ast})$}.
We obtain an isomorphism
$$\overline{\rho}(E_r^{p+r\ast,q-(r-1)\ast}): ~
\otimes_n \det(H^n(E_r^{p+r\ast,q-(r-1)\ast})^{(-1)^n}
\longrightarrow
\otimes_n \det(E_r^{p+rn,q-(r-1)n})^{-1)^n}.$$
Taking \mbox{$\otimes_{p=0}^{r-1} \otimes_{q\ge 0}
\left(\overline{\rho}(E_r^{p+r\ast,q-(r-1)\ast})\right)^{(-1)^{p+q}}$}
yields an isomorphism for $r \ge 2$
$$\overline{\rho}:
\det((H^{0}(E_{r}^{p+r\ast,q-(r-1)\ast}))_{p,q})
\longrightarrow \det((E_{r+1}^{p,q})_{p,q}).$$

Together with the natural isomorphisms
 $V^{p,q}: H^p_{\dR}(B,H^q_{\dR}(F,V)\to E_2^{p,q}$ and $\phi_r^{p,q}$
 between the homology of $E_r$ and $E_{r+1}$, this
definies an isomorphism
\begin{eqnarray}
\overline{\rho}^{\Serre}_{\dR}(f):~
\det((H^p_{\dR}(B;H^q_{\dR}(F;V))_{p,q})
&\xrightarrow{\cong} &
\det((H^n_{dR}(E;V))_{n})
\label{rho^Serre for determinants}
\end{eqnarray}
as the composition of isomorphisms (use the fact that $E_r$
and $E_{\infty}$ are equal for $r$ large enough)
\begin{center}
$\det((H^p_{\dR}(B;H^q_{\dR}(F;V))_{p,q})
\xrightarrow{(V^{p,q})_{p,q}}
\det((E_2^{p,q})_{p,q})
\xrightarrow {\overline{\rho}}
\det((H^{0}(E_2^{p+2\ast,q-\ast}))_{p,q})
\xrightarrow{(\phi_2^{p,q})_{p,q}}
\det((E_3^{p,q})_{p,q})
\xrightarrow{\overline{\rho}}
\det((H^{0}(E_3^{p+3\ast,q-2\ast}))_{p,q})
\xrightarrow{(\phi_3^{p,q})_{p,q}}
\det((E_4^{p,q})_{p,q})
\longrightarrow \ldots \longrightarrow
\det((E_{\infty}^{p,q})_{p,q})
\xrightarrow{((\psi^{p,q})^{-1})_{p,q}}
\det((F^{p,q}/F^{p+1,q-1})_{p,q})
\xrightarrow {\overline{\rho}^{\fil}}
\det((H^n_{\dR}(E;V))_n))$.
\end{center}

For the spectral sequence associated to an arbitrary filtration
$F_*C$ of a finite vector space $W$ one can define 
 $\overline{\rho}_{\fil}(F_*C)$ and $\overline{\rho}^{\ge 2}_{\fil}(F_*C)$ 
in the obvious analogous way.

Now, we have to relate the invariants $\overline{\rho}$ to the previously 
defined real numbers $\rho$ ---using given inner products:.

If the finite-dimensional vector space $V$ has a Hilbert structure,
it induces in a canonical way a Hilbert space structure on $\det(V)$.
In particular the norm $\|v\| \in \rr^{\ge 0}$  is defined for any
element $v \in V$. If $V$ and $W$ come with Hilbert space structures,
for any element
$$u \in \hom_{\rr}(\det(U),\det(V)) = \det(V)^{\ast} \otimes \det(W)$$
we get its norm
\begin{eqnarray}
\|u\| & \in & \rr^{\ge 0}.
\label{norm of an element}
\end{eqnarray}
It turns out that for a finite Hilbert cochain complex $C$ with given Hilbert
structures on $H(C)$
\[ \rho(C)=\ln(\|\overline{\rho}(C)\|) \]

In the case of the fibration $f:E\to B$ of the introduction,
the Riemannian metrics on $E$, $B$ and $V$ induce harmonic Hilbert
structures on $H^n_{\dR}(E;V)$ and $H^p_{\dR}(B;H^q_{\dR}(F;V))$
as described there. Then the norm of the
element $\overline{\rho}^{\Serre}_{\dR}(f)$
defined in \ref{rho^Serre for determinants} is just
\begin{eqnarray}
\rho^{\Serre}_{\dR}(f) & = &
\ln\left(\|\overline{\rho}^{\Serre}_{\dR}(f)\|\right) .
\label{definition of rho^Serre as real number}
\end{eqnarray}

Of course, with similar definitions
a similar result can be obtained for $\rho^{\Serre}_{\sing}(f)$. Also, the
corresponding equation for $\rho_{\fil}$ of a filtration is true.
\end{appendix}


\typeout{-------------------------- references
 --------------------------}

\begin{center}
Address\\
Wolfgang L\"uck und Thomas Schick \\
Fachbereich f\"ur Mathematik und Informatik,
\\
Westf\"alische Wilhelms-Universit\"at M\"unster,
\\
Einsteinstr. 62,
48149 M\"unster,
Germany,
\\
email: lueck@math.uni-muenster.de, schickt@math.uni-muenster.de\\
Fax: 0251 83-38370
\end{center}
\begin{center} Version of \today \end{center}

\end{document}